
%
\parindent=0mm
\font \fiverm=cmr5
\font \sixrm=cmr6
\font \sevenrm=cmr7
\font \eightrm=cmr8
\font \ninerm=cmr9

\font \ninebf=cmbx9
\font \bigbf=cmbx10 scaled \magstep1
\font \Bigbf=cmbx10 scaled \magstep2

\font \tengoth=eufm10
\font \sevengoth=eufm7
\font \fivegoth=eufm5

\newfam\gothfam
\textfont \gothfam=\tengoth
\scriptfont \gothfam=\sevengoth
\scriptscriptfont \gothfam=\fivegoth

%
%
\newfam\srmfam
\textfont \srmfam=\eightrm
\scriptfont \srmfam=\sixrm
\scriptscriptfont \srmfam=\fiverm

\font \tengothb=eufb10
\font \sevengothb=eufb7
\font \fivegothb=eufb5

\newfam\gothbfam
\textfont \gothbfam=\tengothb
\scriptfont \gothbfam=\sevengothb
\scriptscriptfont \gothbfam=\fivegothb

\font \tenmath=msbm10
\font \sevenmath=msbm7
\font \fivemath=msbm5

\newfam\mathfam
\textfont \mathfam=\tenmath
\scriptfont \mathfam=\sevenmath
\scriptscriptfont \mathfam=\fivemath
\def\math{\fam\mathfam\tenmath}
%
%
\def\titre#1{\centerline{\Bigbf #1}\nobreak\nobreak\vglue 10mm\nobreak}

\def\paragraphe#1{\bigskip\goodbreak {\bigbf #1}\nobreak\vglue 12pt\nobreak}
\def\alinea#1{\medskip\allowbreak{\bf#1}\nobreak\vglue 9pt\nobreak}
\def\ssq{\smallskip\qquad}

%
%
\def\th#1{\bigskip\goodbreak {\bf Theor\`eme #1.} \par\nobreak \sl }
\def\prop#1{\bigskip\goodbreak {\bf Proposition #1.} \par\nobreak \sl }
\def\lemme#1{\bigskip\goodbreak {\bf Lemme #1.} \par\nobreak \sl }
\def\cor#1{\bigskip\goodbreak {\bf Corollaire #1.} \par\nobreak \sl }
\def\dem{\bigskip\goodbreak \it D\'emonstration. \rm}
\def\ndem{\bigskip\goodbreak \rm}
\def\qed{\par\nobreak\hfill $\bullet$ \par\goodbreak}
%
%
\def\uple#1#2{#1_1,\ldots ,{#1}_{#2}}
\def\corde#1#2-#3{{#1}_{#2},\ldots ,{#1}_{#3}}

\def\ordcorde#1#2-#3{{#1}_{#2} \le \cdots \le {#1}_{#3}}
\def\strictordcorde#1#2-#3{{#1}_{#2} < \cdots < {#1}_{#3}}
\def \restr#1{\mathstrut_{\bigl |}\raise-8pt\hbox{$\scriptstyle #1$}}
\def \srestr#1{\mathstrut_{\scriptstyle |}\hbox to -1.5pt{}\raise-4pt\hbox{$\scriptscriptstyle #1$}}
\def \inver{^{-1}}
\def\dbar{d\!\!\hbox to 4.5pt{\hfill\vrule height 5.5pt depth -5.3pt
	width 3.5pt}}

\def\frac#1#2{{\textstyle {#1\over #2}}}

\def\R{{\math R}}
\def\C{{\math C}}

\def\N{{\math N}}

\def\fleche#1{\mathop{\hbox to #1 mm{\rightarrowfill}}\limits}
\def\gfleche#1{\mathop{\hbox to #1 mm{\leftarrowfill}}\limits}
\def\inj#1{\mathop{\hbox to #1 mm{$\lhook\joinrel$\rightarrowfill}}\limits}
\def\ginj#1{\mathop{\hbox to #1 mm{\leftarrowfill$\joinrel\rhook$}}\limits}
\def\surj#1{\mathop{\hbox to #1 mm{\rightarrowfill\hskip 2pt\llap{$\rightarrow$}}}\limits}
\def\gsurj#1{\mathop{\hbox to #1 mm{\rlap{$\leftarrow$}\hskip 2pt \leftarrowfill}}\limits}
%
%
\def \g#1{\hbox{\tengoth #1}}
\def \sg#1{\hbox{\sevengoth #1}}

\def\Cal #1{{\cal #1}}
%
%

\def \mop#1{\mathop{\hbox{\rm #1}}\nolimits}
\def \smop#1{\mathop{\hbox{\sevenrm #1}}\nolimits}

\def \mopl#1{\mathop{\hbox{\rm #1}}\limits}
\def \smopl#1{\mathop{\hbox{\sevenrm #1}}\limits}

%
%
\def \bib #1{\null\medskip \strut\llap{[#1]\quad}}
\def\cite#1{[#1]}

\magnification=\magstep1
\parindent=0cm
\def\titre#1{\centerline{\Bigbf #1}\vskip 16pt}
\def\paragraphe#1{\bigskip {\bigbf #1}\vskip 12pt}
\def\alinea#1{\medskip{\bf #1}\vskip 6pt}
\def\ssq{\smallskip\qquad}

\input epsf
\long\def\dessin#1#2{\null
                \bigskip
                \begingroup \epsfysize = #1 $$\epsfbox {#2}$$ \endgroup
                \bigskip
                \goodbreak}                
\titre{Front d'onde et propagation des singularit\'es}
\vskip -4mm
\titre{pour un vecteur-distribution}
\vskip 2mm
\centerline{Dominique MANCHON}
\smallskip
\centerline{Institut Elie Cartan, Universit\'e Henri Poincar\'e - CNRS - INRIA}
\smallskip
\centerline{BP 239, 54506 Vand\oe uvre les Nancy Cedex}
\smallskip
\centerline{\tt manchon@iecn.u-nancy.fr}
\vskip 6mm plus 2mm
\centerline{\bigbf Table des mati\`eres}
\medskip
\+ Introduction  \hbox to 10cm{}                                        &2\cr
\smallskip
\+ I. Calcul fonctionnel holomorphe pour des symboles elliptiques       &5\cr
\+ \qquad I.1. Rappels sur le calcul symbolique  &       5\cr
\+ \qquad I.2. Une r\'esolvante approch\'ee      &       7\cr
\+ \qquad I.3. Calcul fonctionnel holomorphe     &       10\cr
\smallskip
\+ II. Front d'onde d'un vecteur-distribution   & 14\cr
\+ \qquad II.1. Espaces de Sobolev-Goodman       &       14\cr
\+ \qquad II.2. Op\'erateurs pseudo-diff\'erentiels sur les espaces de repr\'esentations                                         & 15\cr
\+ \qquad II.3. D\'efinition et premi\`eres propri\'et\'es du front d'onde & 17\cr
\+ \qquad II.4. Une autre caract\'erisation du front d'onde      & 19\cr
\+ \qquad II.5. Restrictions                                    & 21 \cr
\smallskip
\+ III. Propagation des singularit\'es                          & 23\cr
\+ \qquad III.1. Symboles classiques                            & 23\cr
\+ \qquad III.2. Une forme faible de l'in\'egalit\'e de G\aa rding pr\'ecis\'ee                                 &23 \cr
\+ \qquad III.3. Th\'eor\`eme de propagation des singularit\'es &24\cr
\smallskip
\+ Appendice~: calcul symbolique \`a param\`etre        & 29\cr
\smallskip
\+ R\'ef\'erences                                       & 35\cr
\vskip 6mm minus 2mm
\begingroup\ninerm
{\ninebf R\'esum\'e}~: Nous d\'efinissons le front d'onde d'un vecteur-distribution pour une repr\'esentation unitaire d'un groupe de Lie r\'eel $G$ \`a l'aide du calcul pseudo-diff\'erentiel mis au point dans un travail ant\'erieur [M2]. Cette notion pr\'ecise celle de front d'onde d'une repr\'esentation introduite par R. Howe [Hw]. En application nous donnons une condition suffisante pour qu'un vecteur-distribution reste un vecteur-distribution pour la restriction de la repr\'esentation \`a un sous-groupe ferm\'e $H$, et nous donnons un th\'eor\`eme de propagation des singularit\'es pour les vecteurs-distribution.
\medskip
{\ninebf Abstract}~: We define the wave front set of a distribution vector of a unitary representation in terms of pseudo-differential-like operators [M2]for any real Lie group $G$. This refines the notion of wave front set of a representation introduced by R. Howe [Hw]. We give as an application a necessary condition so that a distribution vector remains a distribution vector for the restriction of the representation to a closed subgroup $H$, and we give a propagation of singularities theorem for distribution vectors.
\endgroup
\smallskip
{\bf Mathematics Subject Classification}~: 22E30, 38S05, 58G15
\paragraphe{Introduction}
\qquad
Cet article est la suite d'un travail r\'ecent \cite {M4} dans lequel \'etait d\'evelopp\'ee la notion de front d'onde d'une repr\'esentation unitaire \cite {Hw}. Nous \'etendons ici la notion de front d'onde aux vecteurs-distribution.
\ssq
La d\'efinition du front d'onde d'une distribution sur une vari\'et\'e $M$ est due \`a L. H\"ormander, qui en donne \'egalement une caract\'erisation \`a l'aide des op\'erateurs pseudo-diff\'erentiels \cite {Hr2 th. 18.1.27}~: le front d'onde est une partie conique ferm\'ee du fibr\'e cotangent priv\'e de la section nulle $T^*M-\{0\}$, donn\'ee par~:
$$WF(u)=\bigcap \mop{char} P,\leqno{(*)}$$
o\`u $P$ parcourt l'ensemble des op\'erateurs pseudo-diff\'erentiels d'ordre z\'ero tels que $Pu\in C^\infty(M)$, et $\mop{char}P$ d\'esigne l'ensemble des $(x,\xi)\in T^*M-\{0\}$ tels que le symbole $p$ de $P$ v\'erifie~: $\mop{lim inf}p(x,t\xi)=0$ lorsque $t\rightarrow +\infty$ (voir aussi \cite {T \S \ VI.1}). Pour un r\'eel $s$ quelconque le front d'onde d'ordre $s$ de la distribution $u$ est d\'efini de la m\^eme mani\`ere, l'intersection s'\'etendant \`a tous les o.p.d. $P$ d'ordre z\'ero tels que $Pu$ appartienne \`a l'espace de Sobolev $H_s^{\smop{loc}}(M)$ \cite{D-H \S\ 6.1}.
\ssq
Or nous avons d\'efini (\cite {M2}, cf. \S\ II.2) pour toute repr\'esentation unitaire $\pi$ d'un groupe de Lie r\'eel $G$ des op\'erateurs pseudo-diff\'erentiels sur l'espace $\Cal H_\pi$ de la repr\'esentation, dont le symbole est une fonction $C^\infty$ sur la vari\'et\'e de Poisson lin\'eaire $\g g^*$ (ces op\'erateurs sont l'image par $\pi$ d'op\'erateurs pseudo-diff\'erentiels invariants sur le groupe, cf. \cite {Me1}, \cite {Stk})~: pour tout $p\in C^\infty(\g g^*)$ tel que sa transform\'ee de Fourier soit une distribution \`a support compact contenu dans un voisinage exponentiel de $0$ dans l'alg\`ebre de Lie $\g g$ on d\'efinit l'op\'erateur $p^{W,\pi}$ de domaine $H_\pi^\infty$ par~:
$$<p^{W,\pi}u,\,v>=<\varphi,C_{u,v>}$$
pour tout $u\in \Cal H_\pi^\infty$ et $v\in \Cal H_\pi$. Ici $C_{u,v}$ d\'esigne le coefficient $<\pi(.)u,\,v>$ et $\varphi$ est donn\'ee par~:
$$j_G.\exp^*\varphi=\Cal F\inver p,$$ 
o\`u $j_G$ d\'esigne le jacobien de l'exponentielle. Il est donc naturel de d\'efinir (\S\ II.3) le front d'onde d'un vecteur-distribution par une formule analogue \`a (*)~:
$$WF(u)=\bigcap \mop{char} p, \leqno{(**)}$$
o\`u $p$ parcourt l'ensemble des symboles d'ordre z\'ero tels que $p^{W,\pi}u$ soit un vecteur $C^\infty$, et o\`u $\mop{char} p$ d\'esigne l'ensemble des directions dans $\g g^*$ o\`u $p$ s'annule \`a l'infini. Nous d\'efinissons \'egalement le front d'onde d'ordre $s$ \`a l'aide des espaces de Sobolev-Goodman $\Cal H_\pi^s$ \cite {Go} dont nous rappelons la d\'efinition au \S\ II.1.
\ssq
Le front d'onde d'un vecteur-distribution est un c\^one ferm\'e $\mop{Ad}^*G$-invariant de $\g g^*-\{0\}$. Ceci raffine la notion de front d'onde d'une repr\'esentation unitaire introduite par R. Howe \cite{Hw}, en ce sens que l'on a l'inclusion~:
$$WF(u)\subset -WF_\pi^e.$$
Nous donnons \'egalement au \S\ II.4 une caract\'erisation du front d'onde de R. Howe en termes d'op\'erateurs pseudo-diff\'erentiels~:
$$WF_\pi^e=-\bigcap \mop{char p},$$
o\`u $p$ parcourt l'ensemble des symboles d'ordre z\'ero tels que $p^{W,\pi}$ soit un op\'erateur r\'egulari\-sant. Le signe moins est un simple artifice de convention.
\ssq
En application nous donnons au \S\ II.5 une condition suffisante pour qu'un vecteur-distribution soit encore un vecteur-distribution pour la restriction de $\pi$ \`a un sous-groupe ferm\'e $H$~:
\th{II.5.2}
1). Soit $u\in \Cal H_{\pi\srestr H}^\infty$. Consid\'erant alors $u$ comme un vecteur-distribution de $\pi$ on a l'inclusion~:
$$WF(u)\subset \g h^\perp.$$
2). Soit $v\in\Cal H_\pi^{-\infty}$. Alors si $WF(v)\cap \g h^\perp=\emptyset$ le vecteur-distribution $v$ est un vecteur-distribution pour la restriction $\pi\restr H$.
\ndem
Nous donnons \'egalement une condition suffisante sur une distribution $\varphi$ \`a support compact sur $G$ pour que $\pi(\varphi)u$ soit un vecteur $C^\infty$~: d\'efinissant $\Cal {WF}(u)$ comme la partie conique de $T^*G-\{0\}$ form\'ee par les translat\'es (\`a gauche ou \`a droite) de $WF(u)$ on a le r\'esultat suivant~:
\prop{II.4.4}
Supposons que $\Cal H_\pi$ soit un espace de Hilbert s\'eparable. Soit $u\in\Cal H_\pi^{-\infty}$. Alors pour toute distribution $\varphi\in\Cal E'(G)$ telle que~:
$$WF(\varphi)\cap \Cal {WF}(u)=\emptyset,$$
on a~:
$$\pi(\varphi)u\in\Cal H_\pi^\infty.$$
\ndem
Le r\'esultat central de cet article est le th\'eor\`eme de propagation des singularit\'es III.3.1, qui entra\^\i ne (corollaire III.3.5) que pour un symbole \`a valeurs r\'eelles $p$ la diff\'erence des deux fronts d'onde~:
$$WF(u)\backslash WF(p^{W,\pi}u)$$
est invariante par le flot du champ hamiltonien d\'efini par le symbole principal $p_m$. La d\'emonstration est parall\`ele \`a celle donn\'ee dans \cite {T\S\ VI.1} pour les distributions, et repose sur une forme faible de l'in\'egalit\'e de G\aa rding pr\'ecis\'ee (Th\'eor\`eme III.2.1).
\ssq
La premi\`ere partie est consacr\'ee au calcul fonctionnel holomorphe pour des symboles elliptiques. Nous avons repris la construction de M.A. Shubin \cite{Shu} en suivant la suggestion de R. Strichartz \cite {St} de consid\'erer des fonctions holomorphes dans un secteur du plan complexe plus g\'en\'erales que $t\mapsto t^z$. Le d\'efaut de cette approche est de limiter le calcul fonctionnel holomorphe \`a des symboles polynomiaux (ou, si l'on veut, aux op\'erateurs diff\'erentiels). La m\'ethode de R. Seeley \cite S s'applique \`a tout o.p.d. elliptique, mais nous n'avons pas su adapter sa m\'ethode \`a notre cadre. La difficult\'e provient de la n\'ecessit\'e de se limiter \`a des symboles dont la transform\'ee de Fourier inverse est \`a support compact, ce qui oblige \`a couper les hautes fr\'equences \`a chaque \'etape de la construction de la r\'esolvante approch\'ee. Cette restriction n'ayant pas lieu d'\^etre dans le cas nilpotent simplement connexe \cite {M3} on peut donc d\'efinir dans ce cas pr\'ecis un calcul fonctionnel holomorphe pour tout symbole elliptique, polynomial ou non.
\ssq
Ce calcul fonctionnel est n\'ecessaire pour montrer la continuit\'e Sobolev de $\Cal H_\pi^s$ dans $\Cal H_\pi^{s-m}$ des op\'erateurs $p^{W,\pi}$ o\`u $p$ est un symbole d'ordre $m$ (proposition II.2.1), contrairement au cas des distributions sur une vari\'et\'e, o\`u ce r\'esultat se montre directement \cite {T. th. II.6.5}. Il repose quant \`a lui sur une variante \`a param\`etre du calcul symbolique (Th\'eor\`eme I.2.4) sur $\g g^*$, que nous d\'eveloppons en appendice. Nous reprenons \cite {M1 th. 4.2} en incluant le param\`etre et en apportant quelques simplifications et corrections (voir la remarque suivant le th\'eor\`eme I.1.2).
\medskip
{\bf Notations}~:
\smallskip
(0.1). On d\'esignera par $G$ un groupe de Lie r\'eel connexe de dimension $n$, d'alg\`ebre de Lie $\g g$ et de dual $\g g^*$. On d\'esigne par $dx$ une mesure de Lebesgue sur $\g g$, et par $dg$ une mesure de Haar \`a gauche sur $G$ normalis\'ee de telle fa\c con que le jacobien $j_G$ de l'exponentielle s'\'ecrive~:
$$j_G(x)=\bigl|\det({1-e^{-\smop{ad}x} \over \mop{ad}x})\bigr|$$
On d\'esigne par $\Delta_G=\det\mop{Ad}g$ la fonction module. Par le choix de la mesure de Haar $dg$ nous identifierons l'espace $C^\infty(G)$ et l'espace des densit\'es $C^\infty$ sur $G$. Par dualit\'e l'espace des distributions sur $G$ s'identifie \`a l'espace des fonctions g\'en\'eralis\'ees sur $G$.
\smallskip
(0.2). On d\'esignera par $T^*G\backslash \{0\}$ le fibr\'e cotangent de $G$ priv\'e de la section nulle. Une partie $C$ de $T^*G\backslash\{0\}$ est {\sl conique\/} si pour tout $(g,\xi)\in C$, $(g,t\xi) \in C$ pour tout r\'eel $t>0$. On notera alors $-C$ l'ensemble des $\{(g,-\xi),\, (g,\xi)\in C\}$.
\smallskip
(0.3). On d\'esignera par $\pi$ une repr\'esentation unitaire fortement continue du groupe de Lie $G$ dans un espace de Hilbert s\'eparable qui sera not\'e $\Cal H_\pi$. On d\'esignera alors par $\Cal H^\infty_\pi$ l'espace des vecteurs ind\'efiniment diff\'erentiables de la repr\'esentation. L'espace $\Cal H_\pi^\infty$ est constitu\'e des vecteurs $u$ tels que pour tout $v\in\Cal H_\pi$ le coefficient~:
$$C_{u,v}:g\longmapsto <\pi(g)u,\,v>$$
soit $C^\infty$ sur $G$. Son dual est l'espace $\Cal H_\pi^{-\infty}$ des vecteurs-distribution.
\smallskip
(0.4). On d\'efinit pour tout $\varphi\in\Cal E'(G)$ l'op\'erateur (en g\'en\'eral non born\'e) $\pi(\varphi)$ de domaine $\Cal H_\pi^\infty$ par la formule~:
$$<\pi(\varphi)u,\, v>=<\varphi,\, C_{u,v}>$$ 
(Voir [J\o] pour une d\'efinition dans le cadre g\'en\'eral des repr\'esentations dans un espace de Banach).
\smallskip
(0.5). Pour toute distribution $\varphi\in \Cal E'(G)$ on pose~:
$$\varphi^*=\Delta_G.\overline {i^*\varphi}$$
o\`u $\Delta_G$ d\'esigne la fonction module et $i$ le diff\'eomorphisme $g\mapsto g\inver$. On a pour tout couple $(u,v)$ dans $\Cal H_\pi^\infty$~:
$$<\pi(\varphi)u,\,v>=<u,\pi(\varphi^*)v>.$$
\paragraphe{I. Calcul fonctionnel holomorphe pour des symboles elliptiques}
\alinea{I.1. Rappels sur le calcul symbolique}
\qquad Soit $G$ un groupe de Lie r\'eel connexe, $\g g$ son alg\`ebre de Lie et $\g g^*$ son dual, que nous verrons comme une vari\'et\'e de Poisson lin\'eaire munie du crochet de Poisson~:
$$\{f,g\}(\xi)=<\xi,\, [Df(\xi),Dg(\xi)]>.$$
Soit $Q$ un voisinage compact de $0$ dans $\g g$. On d\'esigne par $S^m_\rho(\g g^*)$, $m\in\R$, $\rho\in ]0,1]$, la classe de symboles constitu\'ee par les fonctions $p\in C^\infty(\g g^*)$ v\'erifiant les estimations~:
$$|D^\alpha p(\xi)|\le C_\alpha.(1+\|\xi\|^2)^{\frac 12(m-\rho|\alpha|)}.$$
Cet espace, muni des semi-normes $N_\alpha$ d\'efinies comme \'etant les meilleures constantes $C_\alpha$ possibles, est un espace de Fr\'echet. On d\'esignera par $AS^{m,Q}_\rho(\g g^*)$ le sous-espace ferm\'e de $S^m_\rho(\g g^*)$ des symboles $p$ tels que $\mop{supp}\Cal F\inver p\subset Q$, et par $AS^m_\rho(\g g^*)$ la r\'eunion des $AS^{m,Q}_\rho(\g g^*)$, $Q$ parcourant l'ensemble des voisinages compacts de $0$ dans $\g g$. On rappelle deux r\'esultats de \cite {M1} (Theorem 2.1 et Theorem 4.2)~:
\prop{I.1.1 \rm (Th\'eor\`eme d'approximation)}
Soit $Q$ un voisinage compact de $0$ dans $\g g$, et soit $\chi\in C^\infty (\g g)$ \`a support dans $Q$ telle que $\chi=1$ au voisinage de $0$. Alors l'op\'erateur $T^Q$ sur $C^\infty(\g g^*)$ d\'efini par~:
$$T^Q=\Cal F\circ \chi\circ \Cal F\inver$$
est continu de $S^m_\rho(\g g^*)$ dans $AS^{m,Q}_\rho(\g g^*)$, et $I-T^Q$ est continu de $S^m_\rho(\g g^*)$ dans l'espace de Schwartz $S(\g g^*)$.
\ndem
\prop{I.1.2 \rm (Calcul symbolique)}
Soit $K$ un voisinage compact exponentiel de $0$ dans $\g g$, assez petit pour que $(\exp K)^2$ soit l'image diff\'eomorphe d'un voisinage $K^2$ de $0$ dans $\g g$. Alors si $\rho>\frac 12$ il existe un voisinage compact $Q$ de $0$ dans $\g g$ tel que la loi $\#$ d\'efinie \`a partir de la convolution $*$ sur le groupe par~:
$$\exp_*(j_G\inver\Cal F\inver p)\circ \exp_*(j_G\inver\Cal F\inver q)
=\exp_*(j_G\inver\Cal F\inver (p\# q))$$
s'\'etend en une correspondance bilin\'eaire continue~:
$$AS^{m_1,Q}_\rho(\g g^*)\times AS^{m_2,Q}_\rho(\g g^*)
\longrightarrow  AS^{m_1+m_2,Q^2}_\rho(\g g^*).$$
De plus on a le d\'eveloppement asymptotique~:
$$p\#q=\sum_0^N C_k(p,q)+R_N(p,q)$$
o\`u les $C_k$ sont des op\'erateurs bi-diff\'erentiels avec~:
$$C_0(p,q)=pq, \hbox to 8mm{} C_1(p,q)=\frac i2 \{p,q\},$$
les $C_k$ (resp. $R_N$) \'etant par ailleurs continus de $AS^{m_1,Q}_\rho(\g g^*)\times AS^{m_2,Q}_\rho(\g g^*)$ dans la classe $AS^{m_1+m_2-k(2\rho-1),Q+Q}_\rho(\g g^*)$ (resp. $AS^{m_1+m_2-(N+1)(2\rho-1),(Q+Q)\cup Q^2}_\rho(\g g^*)$).
\ndem
{\it Remarque\/}~: le th\'eor\`eme n'est pas vrai a priori pour le voisinage compact $K$ lui-m\^eme, contrairement \`a ce qui est dit par erreur dans l'\'enonc\'e du th\'eor\`eme 4.2 de \cite {M1}. La d\'emonstration du lemme III.8 de \cite {M1} n\'ecessite en effet un voisinage compact $Q$ plus petit (voir le lemme A.6 pour la variante avec param\`etre), ce qui n'est nullement restrictif dans les applications.
\ssq 
On a \'egalement une variante \`a param\`etre du th\'eor\`eme d'approximation, tir\'ee de \cite{M2} (Th\'eor\`eme I.2.2)~:
\prop{I.1.3 \rm(Th\'eor\`eme d'approximation \`a param\`etre)}
Soit $(p_\lambda)_{\lambda\in \Cal P}$ une famille de symboles, o\`u l'ensemble de param\`etres $\Cal P$ est c\^one dans un espace vectoriel norm\'e (typiquement, $\Cal P=\R^+$ ou un secteur de \C). On suppose que $p_\lambda$ appartient \`a $S^m_\rho(\g g^*)$ pour tout $\lambda$, et qu'on a en plus les estimations~:
$$|D^\alpha p_\lambda(\xi)|\le C_\alpha(1+\|\xi\|^2)^{\frac 12(m-\rho|\alpha|)}(1+\|\lambda\|^2)^{\frac 12(m'-\rho'|\alpha|)}$$
avec $m'\in\R$ et $\rho'>0$. Alors pour tout $N\in\N$ et tout multi-indice $\alpha$ il existe $C_N,\alpha$ telle que~:
$$|D^\alpha\bigl((I-T^Q)p_\lambda\bigr)(\xi)
        \le C_{N,\alpha}(1+\|\xi\|^2)^{-N}(1+\|\lambda\|^2)^{-N}.$$
\ndem
Le passage du cas $\Cal P=\R^+$ trait\'e dans \cite {M2} au cas g\'en\'eral ci-dessus est imm\'ediat.
\alinea{I.2. Une r\'esolvante approch\'ee}
\qquad Soit $p\in AS^{m,Q}_\rho(\g g^*)$ avec $m>0$, $\rho>\frac 12$ et $Q$ voisinage compact de $0$ dans $\g g$ satisfaisant les hypoth\`eses de la proposition I.1.2. On suppose que $p$ est elliptique, c'est-\`a-dire qu'il existe $R>0$ tel que pour $\|\xi\|>R$ on a~:
$$C_1(1+\|\xi\|^2)^{\frac m2}\le |p(\xi)|\le C_2(1+\|\xi\|^2)^{\frac m2}.$$
Soit $\Cal P$ un secteur angulaire de $\C$ tel que pour $\|\xi\|$ assez grand $p(\xi)$ ne prenne pas ses valeurs dans un voisinage conique de $\Cal P$. Nous allons construire une param\'etrixe pour le symbole \`a param\`etre~:
$$p_\lambda=p-\lambda,\hbox to 8mm{}\lambda\in\Cal P.$$
On posera~:
$$\Lambda(\xi)=(1+\|\xi\|^2)^{\frac 12},\hbox to 8mm{}
        \Lambda_d(\xi,\lambda)=(1+\|\xi\|^2+|\lambda|^{\frac 2d})^{\frac 12}.$$
\lemme{I.2.1}
On a l'estimation pour $|\xi|$ assez grand~:
$$C_1\Lambda_m(\xi,\lambda)^m\le |p_\lambda(\xi)|\le C_2\Lambda_m(\xi,\lambda)^m.$$
\dem
On a l'encadrement~:
$$C_1\Lambda_m(\xi,\lambda)^m\le \Lambda(\xi)^m+|\lambda|\le C_2\Lambda_m(\xi,\lambda)^m.$$
La majoration vient alors de l'in\'egalit\'e~:
$$|p_\lambda(\xi)|\le C.(\Lambda(\xi)^m+|\lambda|).$$
Pour la minoration, on remarque que les hypoth\`eses sur $\Cal P$ entra\^\i nent l'existence d'un $\varepsilon>0$ tel que $|p(\xi)-\lambda|\ge \varepsilon|p(\xi)|$ et $|p(\xi)-\lambda|\ge \varepsilon|\lambda|$. Il existe donc $C'$ telle que~:
$$|p_\lambda(\xi)|\ge C'(\Lambda(\xi)^m+|\lambda|).$$
\qed
\lemme{I.2.2}
Lorsque $p$ est de plus un symbole polynomial, la famille $(p_\lambda)$ v\'erifie les estimations~:
$$|D^\alpha p_\lambda(\xi)|\le C_\alpha\Lambda_m(\xi ,\lambda)^{m-|\alpha|}.$$
\ndem
Ce lemme, \`a d\'emonstration imm\'ediate, est faux si $p$ n'est pas un polyn\^ome~: il est en effet crucial que les d\'eriv\'ees s'annulent \`a partir d'un certain rang. On r\'esume les deux lemmes pr\'ec\'edents en disant que la famille $(p_\lambda)_{\lambda\in\Cal P}$ appartient \`a la classe de symboles \`a param\`etres~:
$$S^m_{1,m}(\g g^*,\Cal P)$$
dont la d\'efinition proc\`ede directement des estimations du lemme 2, et est elliptique en tant que famille \`a un param\`etre de symboles. En utilisant l'ellipticit\'e de $p_\lambda$ et l'expression explicite des d\'eriv\'ees successives de $p_\lambda\inver$ on obtient classiquement le lemme suivant~:
\lemme{I.2.3}
Pour tout $\lambda\in\Cal P$ tel que $|\lambda|>R$ on a les estimations~:
$$|D^\alpha{1\over p_\lambda}(\xi)|\le C_\alpha.\Lambda_m(\xi,\lambda)^{-m-|\alpha|}$$
\ndem
On se donne plus g\'en\'eralement un r\'eel $d>0$, on pose~:
$$\Lambda_d(\xi ,\lambda)=(1+\|\xi\|^2+|\lambda|^{\frac 2d})^{\frac 12}$$
pour $\xi\in\g g^*$ et $\lambda\in\Cal P$, et on d\'efinit la classe de symboles \`a param\`etre $S^m_{\rho,d}(\g g^*,\Cal P)$ comme l'espace des $p_\lambda$ v\'erifiant les estimations~:
$$|D^\alpha p_\lambda (\xi)|\le C_\alpha \Lambda_d(\xi,\lambda)^{m-\rho|\alpha|}$$
Les semi-normes d\'efinies par les meilleures constantes $C_\alpha$ possibles munissent cet espace d'une structure de Fr\'echet. On construit aussi les variantes analytiques $AS^{m,Q}_{\rho,d}$ et $AS^m_{\rho,d}$ de mani\`ere \'evidente. 
Nous aurons besoin d'un calcul symbolique \`a param\`etre, c'est-\`a-dire d'une version \`a param\`etre du th\'eor\`eme I.1.2. La d\'emonstration, strictement parall\`ele \`a celle du th\'eor\`eme 4.2 de \cite {M1}, est donn\'ee en appendice.  
\th{I.2.4 \rm (calcul symbolique \`a param\`etre)}
Soit $K$ un voisinage compact exponentiel de $0$ dans $\g g$, assez petit pour que $(\exp K)^2$ soit l'image diff\'eomorphe d'un voisinage $K^2$ de $0$ dans $\g g$. Alors si $\rho>\frac 12$ il existe  un voisinage compact $Q$ de $0$ contenu dans $K$ tel que la loi $\#$ s'\'etend en une correspondance bilin\'eaire continue~:
$$AS^{m_1,Q}_{\rho ,d}(\g g^*,\Cal P)\times AS^{m_2,Q}_{\rho ,d}(\g g^*, \Cal P)
\longrightarrow  AS^{m_1+m_2,Q^2}_{\rho, d}(\g g^*, \Cal P).$$
De plus dans le d\'eveloppement asymptotique~:
$$p\#q=\sum_0^N C_k(p,q)+R_N(p,q)$$
les $C_k$ (resp. $R_N$) sont par ailleurs continus de $AS^{m_1,Q}_{\rho, d}(\g g^*, \Cal P)\times AS^{m_2,Q}_{\rho ,d}(\g g^*, \Cal P)$ dans la classe $AS^{m_1+m_2-k(2\rho-1),Q+Q}_{\rho, d}(\g g^*, \Cal P)$ (resp. $AS^{m_1+m_2-(N+1)(2\rho-1),(Q+Q)\cup Q^2}_{\rho, d}(\g g^*, \Cal P)$).
\ndem
\qquad Soit $Q$ un voisinage compact de $0$ dans $\g g$ assez petit (voir \S\ I.1), soit $\chi\in C^\infty(\g g)$ telle que $\chi=1$ au voisinage de $0$ et soit, pour $\lambda\in\Cal P$ et $|\lambda|>R$, le symbole $q^1_\lambda\in AS^{m,Q}_1(\lambda)$ d\'efini par~:
$$q^1_\lambda=T^Q({1\over p_\lambda}).$$
\ndem
\lemme{I.2.5}
Pour tout $N\in\N$ et pour tout multi-indice $\alpha$ on a les estimations~:
$$|D^\alpha (q^1_\lambda-{1\over p_\lambda})|\le C_{\alpha,N}\Lambda_m(\xi,\lambda)^{-N}.$$
\dem
Du lemme I.2.3 on tire imm\'ediatement les estimations~:
$$|D^\alpha {1\over p_\lambda}(\xi)|\le C.\Lambda(\xi)^{-\frac m2 -\frac {|\alpha|}2}\Lambda(\lambda^{\frac 1m})^{-\frac m2 -\frac {|\alpha|}2}.$$
D'apr\`es le th\'eor\`eme d'approximation \`a param\`etre (proposition I.1.3) on a donc les estimations~:
$$\eqalign{|D^\alpha {1\over p_\lambda}(\xi)|   &\le C_{N,\alpha}\Lambda(\xi)^{-N}\Lambda(\lambda^{\frac 1m})^{-N}      \cr
                                                &\le C_{N,\alpha}\Lambda_m(\xi,\lambda)^{-N}.\cr}$$
\qed
On d\'eduit des lemmes I.2.3 et I.2.5 que le symbole \`a param\`etre $q^1_\lambda$ v\'erifie comme $p_\lambda\inver$ les estimations du lemme I.2.3. On supposera, quitte \`a rajouter une constante, que $p_\lambda\not =0$ pour tout $\lambda\in\Cal P$ et pour tout $\lambda$ de module $\ge r$, ce qui permet de construire la r\'esolvante sur le secteur $\Cal P$ en entier ainsi que sur le disque de rayon $r$.
\ssq
\prop{I.2.6}
La famille de symboles~:
$$r_\lambda^1=p_\lambda\# q_\lambda^1-1$$
appartient \`a la classe de symboles \`a param\`etre $AS_{1,m}^{-1,Q}(\g g^*,\Cal P)$
\dem
On \'ecrit~:
$$r_\lambda^1=(p_\lambda\# q_\lambda^1-p_\lambda q_\lambda^1)
        +(p_\lambda q_\lambda^1-1),$$
puis on applique le th\'eor\`eme I.2.4 et le lemme I.2.5.
\qed
Quitte \`a restreidre suffisamment le compact $Q$ on peut construire une param\'etrixe d'ordre $N$ pour $p_\lambda$, en posant~:
$$q_\lambda^N=q_\lambda^1\#\sum_{p=0}^{N-1}(-1)^p(r_\lambda)^{\# p},
\hbox to 12mm{}r_\lambda^N=(r_\lambda^1)^{\# N}.$$
Il est clair que l'on a~:
$$p_\lambda \# q_\lambda^N=1+r_\lambda^N$$
avec $r_\lambda^N\in AS_{1,m}^{-N,Q}$. Un raisonnement standard et un contr\^ole soigneux des supports des transform\'ees de Fourier inverses (\cite {M2} th\'eor\`eme I.3.6) permet de construire une param\'e\-trixe $q_\lambda\in AS^{-m,Q}_{1,m}$ d'ordre infini, c'est-\`a-dire telle que~:
$$r_\lambda = p_\lambda\#q_\lambda -1\in\bigcap_N AS_{1,m}^{-N,Q}(\g g^*,\Cal P).$$
Un petit calcul classique (voir par exemple \cite {He}) montre que la param\'etrixe est bilat\`ere, c'est-\`a-dire que l'on a \'egalement~:
$$r'_\lambda =q_\lambda\#p_\lambda -1\in\bigcap_N AS_{1,m}^{-N,Q}(\g g^*,\Cal P).$$
Le calcul symbolique, l'approximation analytique $T^Q$ et la construction de la param\'etrixe pr\'eservent l'holomorphie par rapport au param\`etre $\lambda$. La r\'esolvante $q_\lambda$ est donc holomorphe en $\lambda$.

\alinea{I.3. Calcul fonctionnel holomorphe}
Soit $p$ un symbole polynomial elliptique de degr\'e $m$ sur $\g g^*$, et soit $\Cal P$ un secteur angulaire de $\C$ tel qu'il existe $R>0$ tel que $p(\xi)\notin \Cal P$ pour $\|\xi\|>R$. Soit $\varphi$ une fonction holomorphe sur un voisinage conique ouvert $\Cal V$ du compl\'ementaire de $\Cal P$ dans $\C -\{0\}$ et v\'erifiant~:
$$|\varphi(z)|\le C.(1+|z|^2)^{\frac s2}$$ 
avec $s<0$. Si $a\notin \Cal P$ et $a>r>0$, et si $s<0$ la formule de Cauchy nous permet d'\'ecrire pour tout $a\in\C-\Cal P$~:
$$\varphi(a)={1\over 2i\pi}\int_\Gamma {\varphi(\lambda)\over \lambda -a}\, d\lambda,$$
\dessin {40 mm}{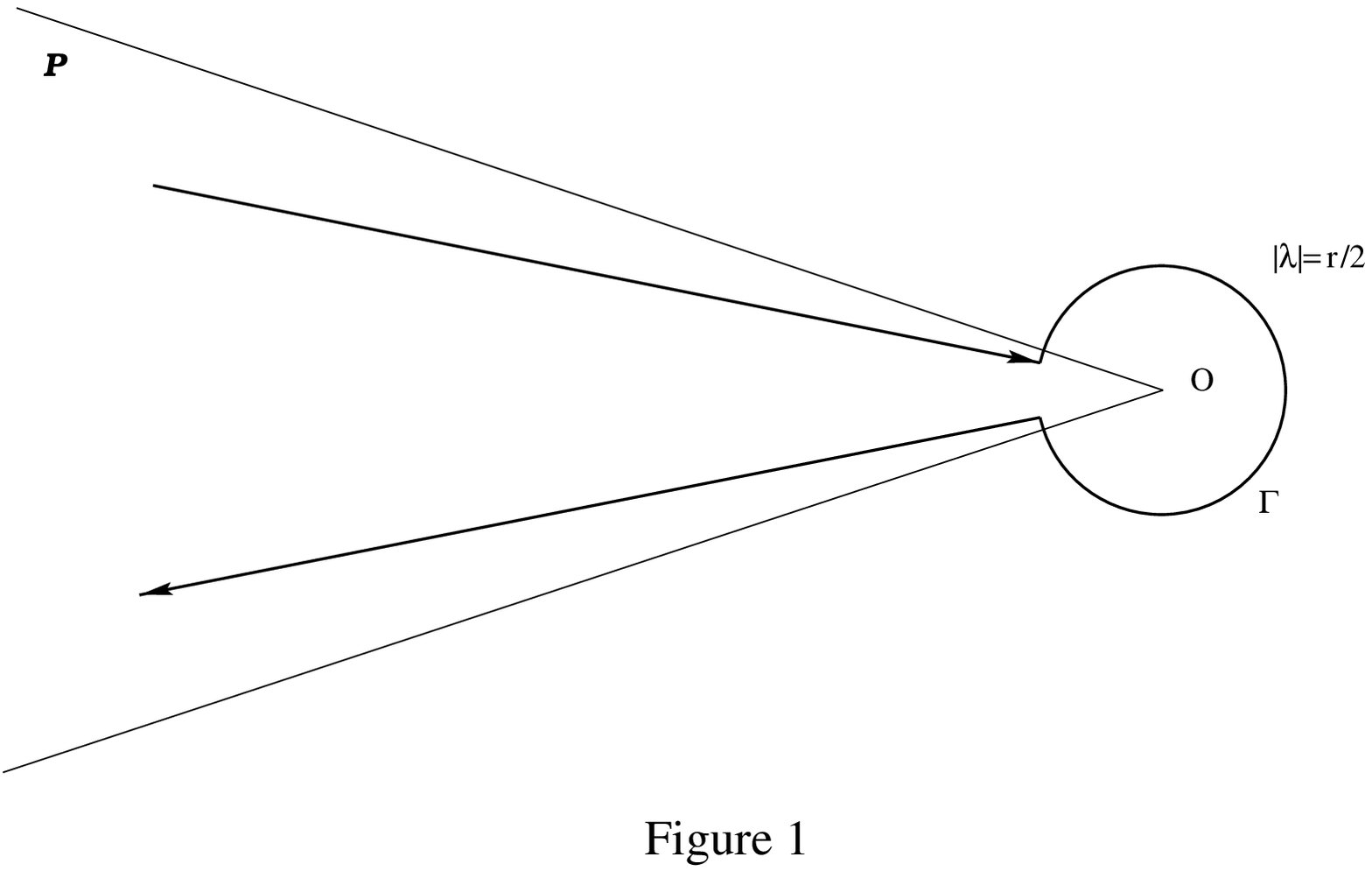}
o\`u $\Gamma\subset \Cal V$ est le contour ci-dessus, d\'efini de la mani\`ere suivante~: quitte \`a faire une rotation on suppose que le secteur $\Cal P$ contient la demi-droite des r\'eels n\'egatifs. On se donne deux r\'eels $-\pi<a<b<\pi$ tels que les deux rayons d'argument $a$ et $b$ soient inclus \`a la fois dans $\Cal V$ et dans l'int\'erieur de $\Cal P$, et on parcourt $\Gamma$ en parcourant d'abord le rayon d'argument $b$ jusqu'au cercle de rayon $\frac r2$ (avec $r<R$), puis le cercle de rayon $\frac r2$ dans le sens n\'egatif, puis le rayon d'argument $a$ en s'\'eloignant de l'origine. Compte tenu de la r\'esolvante approch\'ee construite en I.2 on est donc amen\'e \`a poser~:
\def\de{\mathop{\circ}\limits^\#}
$$\varphi\de p(\xi)=-{1\over 2i\pi}\int_\Gamma \varphi(\lambda)q_\lambda (\xi)
\,d\lambda.$$
Cette formule a un sens, la r\'esolvante $q_\lambda$ construite au \S\ I.2 \'etant d\'efinie pour tout $\lambda\in\Cal P \cup B(0,R)$, donc pour tout $\lambda\in \Gamma$. 
\prop{I.3.1}
$\varphi\de p$ appartient \`a la classe de symboles~:
$$AS^{ms,Q}_1(\g g^*).$$
\dem
dans la formule~:
$$D^\alpha \varphi\de p (t\xi)=-{1\over 2i\pi}\int_{\Gamma}
        \varphi(\lambda)D^\alpha q_\lambda(t\xi)\,d\lambda$$
on suppose $t\ge 1$ et on effectue le changement de variables $\lambda\mapsto t^{-m}\lambda$. Le contour $\Gamma$ est chang\'e en $\Gamma_t=t^{-m}\Gamma$ et il est clair que pour toute fonction holomorphe sur $\Cal V$ son int\'egrale sur $\Gamma_t$ est ind\'ependante de $t$. On a donc~:
$$D^\alpha \varphi\de p(t\xi)=-{1\over 2i\pi}t^m\int_\Gamma
        \varphi(t^m\lambda)D^\alpha q_{t^m\lambda}(t\xi)\,d\lambda.$$
On a donc les estimations~:
$$\eqalign{|D^\alpha\varphi\de p(t\xi)| &\le C.t^m\int_\Gamma
        |t^m\lambda|^s(\|t\xi\|^2+|t^m\lambda|^{\frac 2m})^{-\frac 12(m+|\alpha|)}\, d\lambda \cr
                                        &\le C'.t^{ms-|\alpha|}\int_\Gamma
        |\lambda|^s(\|\xi\|^2+|\lambda|^{\frac 2m})^{-\frac 12(m+|\alpha|)}\, d\lambda \cr}$$
d'o\`u le r\'esultat.
\qed
{\sl Remarque}~: il est a priori surprenant que ce r\'esultat ne fasse pas intervenir les d\'eriv\'ees successives de la fonction holomorphe $\varphi$. Mais les estimations~:
$$|\varphi^{(k)}(z)|\le C_k.(1+|z|^2)^{\frac 12(s-k)}$$
que l'on aurait envie d'imposer sont en fait automatiquement v\'erifi\'ees pour $\xi\in\Cal V$. Pour le voir on part de la formule~:
$$\varphi^{(k)}(tz)={1\over 2i\pi}\int_\Gamma
        {\varphi(\xi)\over (\xi-tz)^{k+1}}\,d\xi$$
et on fait le m\^eme changement de variable que dans la d\'emonstration de la proposition I.3.1. On obtient donc~:
$$\eqalign{|\varphi^{(k)}(tz)|    &={t\over 2\pi}\left |\int_\Gamma
        {\varphi(t\xi)\over (t\xi-tz)^{k+1}}\,d\xi\right |      \cr
                                &\le C.t\int_\Gamma
        {|t\xi|^s\over |t\xi-tz|^{k+1}}\,d\xi           \cr
                                &\le C.t^{s-k}\int_\Gamma
        {|\xi|^s\over |\xi-z|^{k+1}}\,d\xi              .\cr}$$

\prop{I.3.2}
Soient deux fonctions $\varphi$ et $\psi$ holomorphes sur $\Cal P$ et v\'erifiant~:
$$|\varphi(z)|\le C_\alpha (1+|z|^2)^{\frac m2},
\hbox to 8mm{}|\psi(z)|\le C'_\alpha (1+|z|^2)^{\frac {m'}2}.$$
avec $m,m'<0$. Alors la diff\'erence~:
$$(\varphi\de p)\#(\psi\de p) -(\varphi\psi)\de p$$
appartient \`a l'espace de Schwartz $\Cal S(\g g^*)$.
\dem
On consid\`ere un contour $\Gamma'$ tr\`es proche de $\Gamma$, mais situ\'e tout entier du c\^ot\'e oppos\'e \`a $\Cal P$ par rapport \`a $\Gamma$, d\'efini en parcourant d'abord le rayon d'argument $b-\varepsilon$ jusqu'au cercle de rayon $r$, puis le cercle de rayon $r$ dans le sens n\'egatif, puis le rayon d'argument $a+\varepsilon$ en s'\'eloignant de l'origine. On choisit $\varepsilon$ assez petit pour que les deux branches de $\Gamma'$ soient incluses dans $\Cal P$~:
\dessin {40mm}{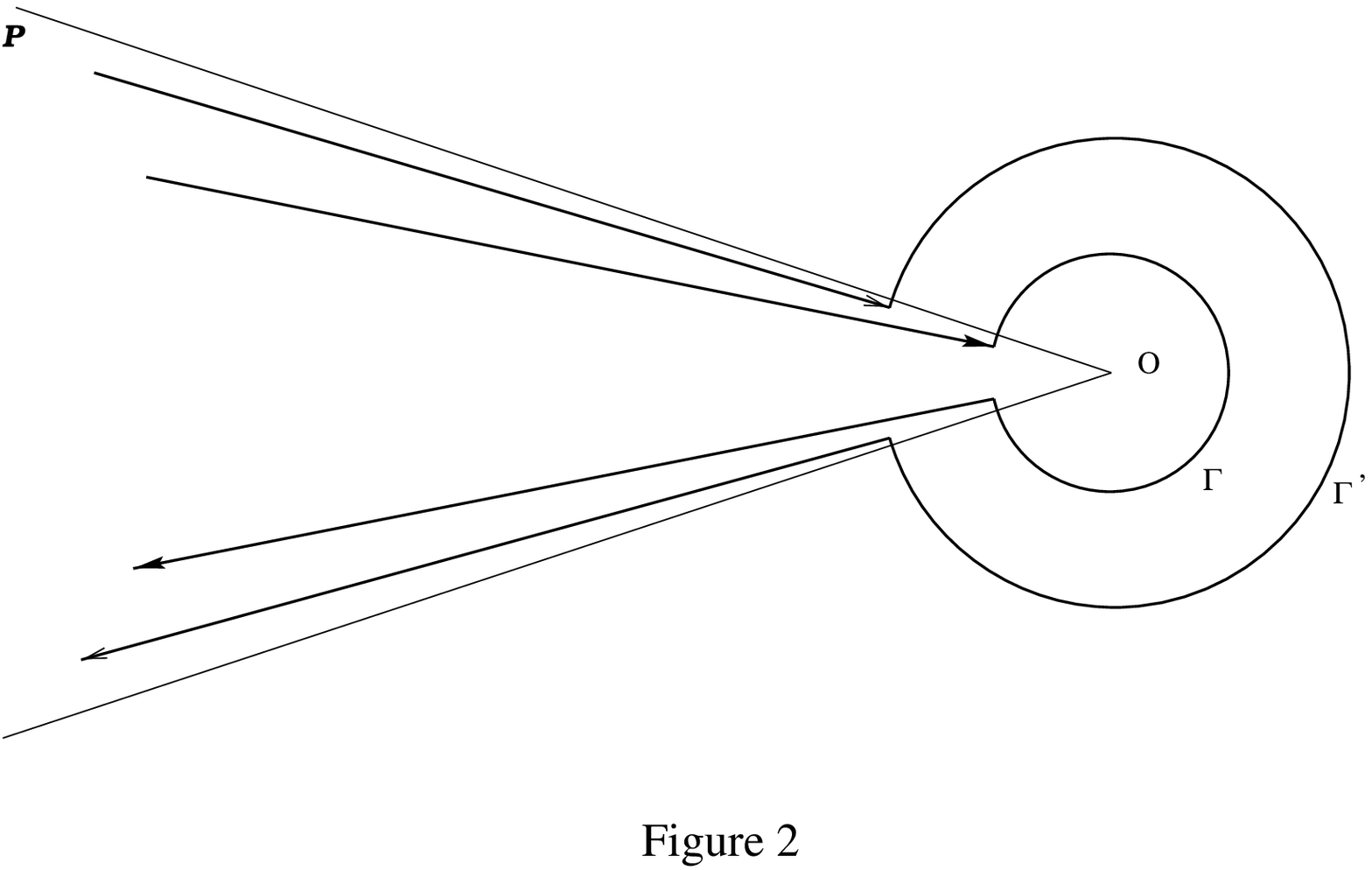}
L'int\'egrale d\'efinissant $\varphi\de p$ ou $\psi\de p$ est la m\^eme en rempla\c cant $\Gamma$ par $\Gamma'$.
\lemme{I.3.3}
Soient $\lambda,\mu\in\C$ avec $\lambda\not =\mu$ et tels que $p$ ne prenne pas ses valeurs en $\mu$ et en $\lambda$. Alors le symbole~:
$$r_{\lambda,\mu}=q_\lambda\#q_\mu-{1\over \lambda -\mu}(q_\lambda-q_\mu)$$
appartient \`a $\Cal S(\g g^*)$, et on a les estimations~:
$$|D^\alpha r_{\lambda,\mu}(\xi)|\le {1\over|\lambda-\mu|}C_{N,\alpha}
        \Lambda (\xi)^{-N},$$
les constantes $C_{N,\alpha}$ \'etant ind\'ependantes de $\lambda$ et $\mu$.
\dem
C'est une cons\'equence imm\'ediate de la formule~:
$$q_\lambda\#q_\mu={1\over \lambda-\mu}q_\lambda \#\bigl((p-\mu)-(p-\lambda)\bigr)\# q_\mu.$$
\qed
On a donc gr\^ace au lemme I.3.3~:
$$\eqalign{\varphi\de p\#\psi\de p      &=-{1\over 4\pi^2}\int_\Gamma\!\int
_{\Gamma'}\varphi(\lambda)\psi(\mu)q_\lambda\#q_\mu\,d\mu\,d\lambda     \cr
                                        &=-{1\over 4\pi^2}
\int_\Gamma\!\int_{\Gamma'}{\varphi(\lambda)\psi(\mu)\over \lambda-\mu}
(q_\lambda-q_\mu)\,d\mu\,d\lambda +\theta,\cr}$$
o\`u $\theta$ appartient \`a l'espace de Schwartz $\Cal S(\g g^*)$. En utilisant la formule de Cauchy pour $\psi(\lambda)$ on obtient donc~:
$$\varphi\de p\#\psi\de p -\theta=
-{1\over 2i\pi}\int_{\Gamma'}\varphi(\lambda)\psi(\lambda)q_\lambda\, d\lambda
        +{1\over 4\pi^2}\int_{\Gamma'}\!\int_{\Gamma}{\varphi(\lambda)
\psi(\mu)\over \lambda-\mu}q_\mu d\lambda \, d\mu.$$
La premi\`ere int\'egrale est \'egale \`a $(\varphi\psi)\de p$ et la deuxi\`eme est nulle. la proposition I.3.2 est donc d\'emontr\'ee.
\qed
\lemme{I.3.4}
si $\varphi(z)=z^{-k}, k\in\N-\{0\}$, le symbole $\varphi\de p -(q_0)^{\#N}$ appartient \`a $\Cal S(\g g^*)$.
\dem
Compte tenu le la proposition I.3.2 il suffit de d\'emontrer le lemme pour $N=1$. On a pour $\varphi(z)=z\inver$:
$$\varphi\de p=-{1\over 2i\pi}\int_{\Gamma''}\lambda\inver q_\lambda\, d\lambda$$
o\`u $\Gamma''$ est un cercle centr\'e \`a l'origine, de rayon petit, et parcouru dans le sens n\'egatif. On a donc $\varphi\de p=q_0$.
\qed
Lorsque l'ordre de croissance $s$ est positif on posera donc~:
$$\varphi\de p =p^{\# k}\# \bigl((z^{-k}\varphi)\de p\bigr).$$
o\`u $k$ est le plus petit entier $\ge -s$. Compte tenu du lemme I.3.4 les propositions I.3.2 et I.3.3 sont encore vraies dans ce cadre. En particulier si on consid\`ere la famille de fonctions $(\varphi_a)_{a\in\C}$ d\'efinie par~:
$$\varphi_a(z)=z^a$$
on a la propri\'et\'e de groupe approch\'ee~:
$$\varphi_{a+b}\de p -(\varphi_a\de p)\# (\varphi_b\de p)\in\Cal S(\g g^*).$$
\smallskip
{\sl Remarque}~: Nous n'avons pas pu d\'efinir $\varphi\de p$ pour un symbole elliptique non polynomial dans $AS_1^{m,Q}$. En effet l'adaptation de la m\'ethode de Seeley \cite {S} n'est pas possible \`a cause de la n\'ecessit\'e de couper les hautes fr\'equences de la param\'etrixe $q^1_\lambda$ d'ordre $1$, ce qui annule le gain en d\'ecroissance par rapport au param\`etre que l'on aurait pu obtenir si on avait pu prendre $q_\lambda^1=(p-\lambda)\inver$ et non une approximation analytique. Or cette op\'eration est inutile dans le cas d'une alg\`ebre de Lie nilpotente (on dispose d'un calcul symbolique avec les classes de symboles ``ordinaires'' $S^m_\rho(\g g^*),\rho>\frac 12$ \cite {M3}). On peut donc dans ce cas suivre \cite S et d\'efinir $\varphi\de p$ pour un symbole elliptique $p\in S^m_\rho(\g g^*)$ quelconque.
\paragraphe{II. Front d'onde d'un vecteur-distribution}
\qquad Nous mettons en \'evidence `a l'aide du calcul holomorphe de la premi\`ere partie la r\'egularit\'e des op\'erateurs pseudo-diff\'erentiels $p^{W,\pi}$ par rapport aux espaces de Sobolev $\Cal H_\pi^s$, dont la d\'efinition est rappel\'ee en II.1. Nous d\'efinissons alors le front d'onde d'un vecteur-distribution, et nous dressons une liste de propri\'et\'es tout-\`a-fait analogues \`a celles du front d'onde d'une distribution. Les techniques employ\'ees sont standard (\cite {D-H} \S\ 6.1, \cite {T} \S\ VI.1). 
\ssq
Nous \'etablissons ensuite le lien avec le front d'onde de la repr\'esentation, et enfin nous appliquons cette construction \`a l'\'etude des vecteurs $C^\infty$ et des vecteurs-distribution de la restriction de la repr\'esentation unitaire \`a un sous-groupe ferm\'e.
\alinea{II.1. Espaces de Sobolev-Goodman}
\qquad Soit $G$ un groupe de Lie r\'eel connexe, $\g g$ son alg\`ebre de Lie et $\g g^*$ son dual. Soit $\pi$ une repr\'esentation unitaire de $G$ dans un espace de Hilbert $\Cal H_\pi$. Soit $\Cal H_\pi^\infty$ l'espace des vecteurs $C^\infty$ de la repr\'esentation, dont le dual $\Cal H_\pi^{-\infty}$ constitue l'espace des vecteurs-distribution.
\ssq
On se donne une base $(X_i)_{i=1,\ldots ,n}$ de $\g g$, et on consid\`ere le laplacien positif~:
$$\Delta =-\sum_1^n X_i^2$$
dans l'alg\`ebre enveloppante $\Cal U(\g g)$. L'op\'erateur non born\'e $\pi(1+\Delta)$ est positif et essentiellement auto-adjoint \cite {Ne}, ce qui permet de d\'efinir ses puissances fractionnaires gr\^ace au th\'eor\`eme spectral.
\ssq
On d\'efinit alors l'{\sl espace de Sobolev-Goodman\/} $\Cal H_\pi^s$ comme le domaine de l'op\'erateur $\pi(1+\Delta)^{\frac s2}$ \cite {Go}. L'espace $\Cal H_\pi^s$, ind\'ependant du choix de la base, est un espace de Hilbert, le produit scalaire \'etant donn\'e par~:
$$<u,\,v>_s=<\pi(1+\Delta)^{\frac s2}u,\,\pi(1+\Delta)^{\frac s2}v>_{\Cal H_\pi}.$$
Dans le cas o\`u $G=\R^n$ et o\`u $\pi$ est sa repr\'esentation r\'eguli\`ere, ces espaces co\"\i ncident avec les espaces de Sobolev usuels $\Cal H_s(\R^n)$. Les $\Cal H_\pi^s$ v\'erifient les propri\'et\'es suivantes~:
\smallskip
1. Si $s\ge t$ alors $\Cal H_\pi^s\subset \Cal H_\pi^t$ et l'inclusion est continue.
\smallskip
2. Le produit scalaire de $\Cal H_\pi$ s'\'etend en une forme sesquilin\'eaire continue~: $\Cal H_\pi^{-s}\times \Cal H_\pi^s \longrightarrow \C$, et de ce fait $\Cal H_\pi^{-s}$ s'identifie au dual de $\Cal H_\pi^s$.
\smallskip
3. L'espace des vecteurs-distribution $\Cal H_\pi^{-\infty}$ est la limite inductive des $\Cal H_\pi^s$ lorsque $s\rightarrow -\infty$, et $\Cal H_\pi^\infty$ est la limite projective des $\Cal H_\pi^s$ lorsque $s\rightarrow +\infty$.
\alinea{II.2. Op\'erateurs pseudo-diff\'erentiels sur les espaces de repr\'esentations}
\qquad On reprend les notations du \S \ I. On associe \cite {M2} \`a tout symbole $p\in AS_\rho^m(\g g^*)$ l'op\'erateur de symbole de Weyl $p$ dans l'espace de la repr\'esentation $\pi$, d\'efini sur $\Cal H_\pi^\infty$ par~:
$$p^{W,\pi}u=\int_{\sg g}\Cal F\inver p(x)\pi(\exp x)\,dx. $$
Cet op\'erateur est en g\'en\'eral non born\'e lorsque $m>0$, mais r\'egularisant (c'est-\`a-dire qu'il transforme tout vecteur en vecteur $C^\infty$) lorsque $p$ appartient \`a l'espace de Schwartz $\Cal S(\g g^*)$. On rappelle \cite {M2} que lorsque $p\in AS^{m_1,Q}_\rho$ et $q\in AS^{m_2,Q}_\rho$ avec $\rho>\frac 12$ et $Q$ assez petit, le produit $p\#q$ est bien d\'efini et on a~:
$$p^{W,\pi}\circ q^{W,\pi}=(p\#q)^{W,\pi}.$$
\prop{II.2.1}
Soit $p\in AS^{m,Q}_\rho(\g g^*)$ avec $\rho>\frac 12$ et $Q$ assez petit. Alors $p^{W,\pi}$ s'\'etend en un op\'erateur continu de $\Cal H_\pi^s$ dans $\Cal H_\pi^{s-m}$.
\dem
Soit $L(\xi)=1+\sum_1^n<\xi,\,X_j>^2$. C'est un symbole polynomial elliptique d'ordre 2, et l'op\'erateur $L^{W,\pi}$ est \'egal \`a $\pi(1+\Delta)$. Soit $q_\lambda$ une r\'esolvante approch\'ee de $L$. On applique le calcul fonctionnel du \S\ I.3 \`a $L$, de sorte que si l'on pose~:
$$p_s=\varphi\de p\hbox{ avec }\varphi(z)=z^{\frac s2}$$
l'op\'erateur $L_s^{W,\pi}-(L^{W,\pi})^{\frac s2}$ est r\'egularisant. On voit ceci en comparant les deux expressions~:
$$L_s^{W,\pi}=-{1\over 2i\pi}\int_\Gamma \lambda^{\frac s2}q_\lambda^{W,\pi}\,d\lambda$$
et
$$(L^{W,\pi})^{\frac s2}=-{1\over 2i\pi}
\int_\Gamma \lambda^{\frac s2} (L^{W,\pi}-\lambda I)\inver\,d\lambda,$$
la deuxi\`eme expression ayant un sens si le rayon $r$ choisi pour la construction du contour $\Gamma$ est assez petit. En particulier $L_s^{W,\pi}$ est continu de $\Cal H_\pi^t$ dans $\Cal H_\pi^{t-s}$ pour tout $t$.
\ssq
\lemme{II.2.2}
Pour tout symbole $p$ d'ordre $0$ et pour tout r\'eel $s$ l'op\'erateur $p^{W,\pi}$ est continu de $\Cal H_\pi^s$ dans $\Cal H_\pi^s$.
\dem
On commence par montrer le lemme pour $0\le s\le 1$. Soit $u\in\Cal H_\pi^\infty$. Les op\'erateurs $p^{W,\pi}$ et $(L_s\#p-p\#L_s)^{W,\pi}$ \'etant continus de $\Cal H_\pi^0$ dans $\Cal H_\pi^0$ on a les estimations suivantes~:
$$\eqalign{\|p^{W,\pi}u\|_s     &\le C.\|L_s^{W,\pi}p^{W,\pi}u\|_0      \cr
        &\le C'.\bigl(\|p^{W,\pi}L_s^{W,\pi}u\|_0+\|(L_s\#p - p\#L_s)^{W,\pi}
u\|_0\bigr)\cr
        &\le C''\|L_s^{W,\pi}u\|_0+C'''\|u\|_0 \cr
        &\le C''''\|u\|_s.\cr}$$
R\'eit\'erant l'argument pr\'ec\'edent en rempla\c cant $\Cal H_\pi^0$ par $\Cal H_\pi^1$ on d\'emontre le lemme pour $1\le s\le 2$, puis de proche en proche pour tout $s\ge 0$. Un simple argument de dualit\'e permet alors d'\'etendre le r\'esultat pour $s\le 0$.
\qed
\qquad On finit la d\'emonstration de la proposition en remarquant que si $p$ est d'ordre $m$ le symbole $L_{-m}\# p$ est d'ordre $0$, ce qui implique \cite {M2 corollaire I.4.2} que $L_{-m}^{W,\pi}\circ p^{W,\pi}$ est born\'e de $\Cal H_\pi$ dans $\Cal H_\pi$. D'apr\`es le lemme II.2.2 cet op\'erateur $L_{-m}^{W,\pi}\circ p^{W,\pi}$ est born\'e de $\Cal H_\pi^s$ dans $\Cal H_\pi^s$ pour tout $s$, donc $L_m^{W,\pi}L_{-m}^{W,\pi}p^{W,\pi}$ est continu de $\Cal H_\pi^s$ dans $\Cal H_\pi^{s-m}$. Comme l'op\'erateur $L_m^{W,\pi}L_{-m}^{W,\pi}$ est
\'egal \`a l'identit\'e modulo un op\'erateur r\'egularisant, ceci termine la d\'emonstration de la proposition II.2.1.
\qed
{\sl Remarque}~: Nous avons exprim\'e les puissances fractionnaires de $L^{W,\pi}$ \`a l'aide d'une int\'egrale sur le contour $\Gamma$, ce qui est possible parce que $L^{W,\pi}$ est auto-adjoint et  $L^{W,\pi}-I$ est positif. Il est possible de faire de m\^eme en rempla\c cant $L$ par n'importe quel symbole $p$ elliptique v\'erifiant les conditions du \S\ I.3, \`a condition que la r\'esolvante soit compacte (ce qui est le cas si la repr\'esentation est fortement tra\c cable, voir par exemple \cite {Hw}, \cite {M4}). En effet dans ce cas le spectre de $p^{W,\pi}$ est discret, son intersection avec $\Cal P$ est born\'ee donc finie, il est donc possible de choisir un contour $\Gamma$ qui \'evite le spectre (voir \cite {Shu \S\ 9}). On d\'efinit ainsi $\varphi(p^{W,\pi})$ pour toute fonction $\varphi$ holomorphe dans le secteur $\Cal P$, et en refaisant le calcul de la proposition I.3.2 on montre que pour deux fonctions holomorphes $\varphi$ et $\psi$ dans $\Cal P$ on a~:
$$\varphi(p^{W,\pi})\circ \psi(p^{W,\pi})=(\varphi\psi)(p^{W,\pi}).$$
\goodbreak
\alinea{II.3. D\'efinition et premi\`eres propri\'et\'es du front d'onde}
\qquad Soit $p\in AS_\rho^0(\g g^*)$ un symbole d'ordre z\'ero. On appelle {\sl ensemble caract\'eristique\/} de $p$ le c\^one d\'efini par~:
$$\mop{char} p=\{\xi\in\g g^*/\mopl{lim inf}_{t\rightarrow +\infty}p(t\xi)=0\}.$$
On montre facilement que si $\rho=1$ ce c\^one est ferm\'e.
On appelle {\sl front d'onde\/} d'un vecteur-distribution $u$ le c\^one ferm\'e de $\g g^*$ d\'efini par~:
$$WF(u)=\bigcap_{{p\in AS_1^0(\g g^*)\atop p^{W,\pi}u\in \Cal H_\pi^\infty}}
\mop{char}p.$$
Le {\sl front d'onde d'ordre $s$\/} du vecteur-distribution $u$ est quant \`a lui d\'efini pour tout r\'eel $s$ par~:
$$WF_s(u)=\bigcap_{{p\in AS_1^0(\g g^*)\atop p^{W,\pi}u\in \Cal H_\pi^s}}
\mop{char}p.$$
Les $WF_s(u)$ forment une famille croissante de c\^ones ferm\'es dont la r\'eunion est le front d'onde  $WF(u)$. Par ailleurs pour tout $g\in G$ on a la relation de covariance~:
$$\pi(g)\circ p^{W,\pi}\circ \pi(g)\inver=(\mop{Ad}^*g.p)^{W,\pi},$$
qui implique que $WF(u)$ et $WF_s(u)$ sont covariants sous l'action coadjointe de $G$ dans $\g g^*$~: pour tout $g\in G$ on a~:
$$WF(\pi(g)u)=\mop{Ad}^*g.WF(u)$$
et de m\^eme pour le front d'onde d'ordre $s$.
\ssq
Le front d'onde d\'efini ci-dessus v\'erifie quelques propri\'et\'es tout \`a fait similaires \`a celles v\'erifi\'ees par le front d'onde d'une distribution. Rappelons \cite{T \S\ VI.1} qu'un symbole $p\in S^m_1(\g g^*)$ est {\sl d'ordre $-\infty$\/} sur un c\^one ouvert $U$ si pour tout c\^one ferm\'e $K$ contenu dans $U$, pour tout $N\in\N$ et pour tout multi-indice $\alpha$ il existe une constante $C_{N,\alpha,K}$ telle que pour tout $\xi\in K$~:
$$|D^\alpha_p(\xi)|\le C_{N,\alpha,K}\Lambda(\xi)^{-N}.$$
Le {\sl support essentiel\/} $ES(p)$ d'un symbole $p\in S^m_1(\g g^*)$ est le plus petit c\^one ferm\'e de $\g g^*-\{0\}$ tel que $p$ soit d'ordre $-\infty$ sur son compl\'ementaire. Le d\'eveloppement asymptotique du produit $\#$ en op\'erateurs bi-diff\'erentiels entra\^\i ne l'inclusion~:
$$ES(p\#q)\subset ES(p)\cap ES(q).$$
Les deux propositions suivantes r\'esument les propri\'et\'es essentielles du front d'onde. La d\'emonstra\-tion, laiss\'ee au lecteur, est strictement parall\`ele \`a celle donn\'ee dans \cite {T \S\ VI.1} dans le cas du front d'onde d'une distribution. L'ingr\'edient essentiel est l'existence d'une param\'etrixe pour les symboles elliptiques.
\ssq
\prop{II.3.1} 
Pour tout c\^one ferm\'e $\Cal C$ tel que $WF(u)\cap \Cal C=\emptyset$ 
il existe un symbole $q=1+a$ avec $ES(a)\cap \Cal C=\emptyset$, tel que~:
$$q^{W,\pi}u\in \Cal H_\pi^\infty.$$
\ndem
\prop{II.3.2}
Soit $u\in \Cal H_\pi^{-\infty}$ et soit $p$ un symbole dans $AS^{m,Q}_1(\g g^*)$, avec $Q$ assez petit. Alors~:
\smallskip
1) Si $ES(p)\cap WF(u)=\emptyset$ alors $p^{W,\pi}u\in\Cal H_\pi^\infty$.
\smallskip
2) $WF(p^{W,\pi}u)\subset WF(u)\cap ES(p)$.
\smallskip
3) $WF(u)\subset WF(p^{W,\pi}u)\cup\mop{char}p$.
\smallskip
4) Si $p$ est elliptique, $WF(p^{W,\pi}u)=WF(u)$.
\ndem
Le front d'onde d'ordre $s$ admet \'egalement la caract\'erisation suivante~:
\prop{II.3.3}
Soit $\xi\in\g g^*-\{0\}$. Alors $\xi\notin WF_s(u)$ si et seulement si on peut \'ecrire $u=u_1+u_2$ avec $u_1\in\Cal H_\pi^s$ et $\xi\notin WF(u)$.
\dem
Si $u$ admet la d\'ecomposition ci-dessus il est clair que $\xi\notin WF_s(u)$. R\'eciproquement si  $\xi\notin WF_s(u)$ il existe un symbole $p$ d'ordre z\'ero tel que $p^{W,\pi}u\in\Cal H^s_\pi$ et $\xi\notin\mop{char}p$. Il existe donc un voisinage conique ferm\'e $\Cal C$ de $\xi$ tel que $p$ soit $\Cal C$-elliptique d'ordre z\'ero, c'est-\`a-dire v\'erifiant~:
$$C_1\le |p(\eta)|\le C_2$$
pour $\eta\in\Cal C$ et $\|\eta\|$ assez grand. La construction de la param\'etrixe (voir \cite {He}, \cite{Hr1}, \cite {M2 prop. I.3.6} et aussi la d\'emonstration de la proposition I.2.6) se restreint sans difficult\'e au c\^one $\Cal C$~: on part d'un symbole $q_1$ d'ordre z\'ero co\"\i ncidant avec $1/p$ sur $\Cal C\cap \{\|\eta\|\ge R\}$, convenablement amput\'e de ses hautes fr\'equences. La construction de \cite {M2 prop. I.3.6} (que nous reprenons dans la d\'emonstration de la proposition I.2.6 dans le cas d'un symbole avec param\`etre) aboutit \`a l'existence d'un symbole $q$ d'ordre z\'ero tel que~:
$$q\# p=1+r$$
avec $ES(r)\cap \Cal C'=\emptyset$ pour tout c\^one ouvert $\Cal C'\subset \Cal C$. On pose alors~:
$$u_1=q^{W,\pi}p^{W,\pi}u\hbox{ et }u_2=r^{W,\pi}u.$$
D'apr\`es la proposition II.2.1 $u_1$ appartient \`a $\Cal H^s_\pi$, et d'apr\`es la proposition II.3.2 $\xi$ n'appar\-tient pas \`a $WF(u_2)$.
\qed
\break
\alinea{II.4. Une autre caract\'erisation du front d'onde}
\prop{II.4.1}
Soit $u\in\Cal H_\pi^\infty$, et soit $T_u$ la distribution sur $G$
\`a valeurs dans $\Cal H_\pi$ d\'efinie par~:
$$T_u=\pi(.)u.$$
Alors le front d'onde de $u$ est \'egal \`a l'oppos\'e du front d'onde en
l'\'el\'ement neutre de la distribution $T_u$.
\dem                     
Les vecteurs $C^\infty$ sont caract\'eris\'es de la fa\c con suivante~:
\lemme{II.4.2}
Pour tout $p\in AS^Q(\g g^*)$, pour tout
$u\in\Cal H_\pi^\infty$ et pour tout $N\in\N$ on a~:
$$\|p_\eta^{W,\pi}u\|\le C_N\|\eta\|^{-N},$$
avec $p_\eta=p(\eta+.)$.
\dem
Si $u$ appartient \`a $\Cal H_\pi^\infty$ la fonction $\psi : x\mapsto \Cal F\inver p(x)\pi(\exp x)u$ est $C^\infty$ \`a support compact, donc sa transform\'ee de Fourier $\widehat\psi$ est \`a d\'ecroissance rapide sur $\g g^*$. Or un petit calcul facile donne~:
$$\widehat \psi(\xi)=p_\xi^{W,\pi}u.$$
\qed
{\sl Fin de la d\'emonstration de la proposition II.4.1\/}~: D'apr\`es la proposition II.3.1 un 
$\xi\in\g g^*-\{0\}$ est hors de $WF(u)$
si et seulement s'il existe un voisinage conique $V_\xi$ de $\xi$ et
un symbole $q=1+a$ avec $ES(a)\cap V_\xi=\emptyset$, tels que~:
$$q^{W,\pi}u\in \Cal H_\pi^\infty.$$
On peut supposer que $V_\xi$ est convexe. Soit maintenant 
$\Cal P$ un c\^one ferm\'e contenu dans $V_\xi$. Pour
tout $\eta\in \Cal P$ le symbole $q_\eta=1+a_\eta$ v\'erifie $q_\eta^{W,\pi}u\in \Cal H_\pi^\infty$, et le support essentiel du symbole \`a param\`etre
$a_\eta$ ne rencontre pas $V_\xi$. Le support essentiel de $r_{-\eta}$ est quant \`a lui contenu dans $V_\xi$. Les supports des transform\'ees de
Fourier inverses de $r_{-\eta}$ et $a_\eta$ sont ind\'ependants de
$\eta$. Supposant ces supports assez petits on peut consid\'erer le
produit $r_{-\eta}\# a_\eta$. Consid\'erant le d\'eveloppement du produit $\#$ en op\'erateurs bi-diff\'erentiels ainsi que l'expression explicite du reste donn\'ee en appendice on montre~:
$$r_{-\eta}\# a_\eta\in AS^{-\infty}(\g g^*,\Cal P).$$
On en d\'eduit que $r_{-\eta}^{W,\pi}a_\eta^{W,\pi} u$ est \`a d\'ecroissance rapide en $\eta$, et donc finalement que $r_{-\eta}^{W,\pi}u$ est aussi \`a d\'ecroissance rapide en $\eta\in\Cal P$, et ce pour tout $r\in AS(\g g^*)$, ce qui veut dire que $\xi$ n'appartient pas \`a $-WF(T_u)$. L'inclusion r\'eciproque, laiss\'ee au lecteur, se d\'emontre en remontant le raisonnement.
\qed
\qquad Il nous reste \`a comparer le front d'onde d'un vecteur-distribution avec le front d'onde d'une repr\'esentation tel qu'il a \'et\'e d\'efini par R. Howe \cite{Hw}. Le front d'onde $WF_\pi$ de la repr\'esentation $\pi$ est une partie conique ferm\'ee de $T^*G-\{0\}$ invariante par translation \`a gauche et \`a droite. Il est donc enti\`erement d\'etermin\'e par son intersection $WF^e_\pi$ avec l'espace cotangent en l'\'el\'ement neutre, qui est un c\^one ferm\'e $\mop{Ad}^*G$-invariant de $\g g^*$. Le crit\`ere {\sl vii}) du th\'eor\`eme 1.4 de \cite{Hw} montre que $WF^e_\pi$ est le front d'onde en $e$ de la distribution $\pi(.)$ \`a valeurs op\'erateurs. L'analogue de la proposition II.4.1, dont la d\'emonstration est identique, nous donne une d\'efinition alternative pour le front d'onde d'une repr\'esentation unitaire~:
\prop{II.4.3}
$$WF_\pi^e=-\bigcap_{{p\in AS_1^0(\g g^*) \atop p^{W,\pi}\hbox{\sevenrm r\'egularisant}}}\mop{char}p.$$
\ndem
En corollaire imm\'ediat nous avons l'inclusion~:
$$WF(u)\subset -WF_\pi^e$$
pour tout vecteur-distribution $u$, et plus pr\'ecis\'ement~:
$$WF_\pi^e=-\bigcup_{u\in\Cal H_\pi^{-\infty}}WF(u).$$
\qquad Soit maintenant $u\in\Cal H_\pi^{-\infty}$, soit $\varphi$ une distribution \`a support compact sur $G$, et soit $\pi(\varphi)$ l'op\'erateur (non born\'e) qu'elle d\'efinit sur l'espace $\Cal H_\pi$. Le domaine de $\pi(\varphi)$ contient $\Cal H_\pi^\infty$. A quelle condition le vecteur $\pi(\varphi)u$ est-il un vecteur $C^\infty$?
\ssq
La r\'eponse \`a cette question est un analogue du th\'eor\`eme III.7 de \cite {M4} pour un vecteur-distribution~: d\'esignons par $\Cal {WF}(u)$ la plus petite partie conique ferm\'ee de $T^*G-\{0\}$ invariante par translation \`a gauche et \`a droite contenant $WF(u)$. 
\prop{II.4.4}
Supposons que $\Cal H_\pi$ soit un espace de Hilbert s\'eparable. Soit $u\in\Cal H_\pi^{-\infty}$. Alors pour toute distribution $\varphi\in\Cal E'(G)$ telle que~:
$$WF(\varphi)\cap \Cal {WF}(u)=\emptyset,$$
on a~:
$$\pi(\varphi)u\in\Cal H_\pi^\infty.$$
\dem
comme la convolution par un \'el\'ement de l'alg\`ebre enveloppante respecte le front d'onde, il suffit de montrer que $\pi(\varphi)u$ appartient \`a $\Cal H_\pi$ pour toute distribution $\varphi\in\Cal E'(G)$ v\'erifiant les conditions du th\'eor\`eme. On proc\`ede alors comme pour le th\'eor\`eme III.4 de \cite {M4}. On remarque d'abord (cf \cite {M4}, prop. III.2 et corollaire III.3) que le front d'onde de la convol\'ee $\varphi^* * \varphi$ ne rencontre pas $\Cal {WF}(u)$. Pour tout $v\in\Cal H_\pi^{-\infty}$ le front d'onde du coefficient-distribution~:
$$C_{u,v}=<\pi(.)u,\,v>$$
est contenu dans $-\Cal {WF}(u)$. C'est en particulier vrai pour le coefficient-distribution diagonal $C_{u,u}$.
\ssq
D\'esignons par $p_t$ le noyau de la chaleur sur le groupe $G$. Lorsque $t$ tend vers $0$ la quantit\'e $\|\pi(p_t)\pi(\varphi)u\|_2^2$ converge en croissant vers $\ell =<(\varphi^**\varphi).C_{u,u},\,1>$, la condition sur les fronts d'onde rendant licite le produit des deux distributions. Par ailleurs si $(e_i)$ d\'esigne une base orthonorm\'ee de $\Cal H_\pi$ form\'ee de vecteurs $C^\infty$, on a~:
$$\|\pi(p_t)\pi(\varphi)u\|^2_2       
                           =\sum_i |<u,\pi(p_t)\pi(\varphi^*)e_i>|^2.$$
Il r\'esulte du th\'eor\`eme de convergence monotone que la somme~:
$$\sum_i|<u,\pi(\varphi^*)e_i>|^2=\sum_i|<\pi(\varphi)u,\, ei>|^2$$
converge vers la limite finie $\ell$ ci-dessus, ce qui veut dire que $u$ appartient \`a $\Cal H_\pi$.
\qed    
\alinea{II.5. Restrictions}
\qquad Soit $G$ un groupe de Lie connexe, $\pi$ une repr\'esentation unitaire de $G$ et $H$ un sous-groupe ferm\'e de $G$. On note comme pr\'ec\'edemment $\Cal H_\pi^\infty$, $\Cal H_\pi^{-\infty}$ et $\Cal H_\pi^s$ l'espace des vecteurs $C^\infty$, des vecteurs-distribution et de Sobolev-Goodman respectivement, pour la repr\'esenta\-tion $\pi$. On note $\Cal H_{\pi\srestr H}^\infty$, $\Cal H_{\pi\srestr H}^{-\infty}$ et $\Cal H_{\pi\srestr H}^s$ les espaces analogues relativement \`a la restriction de $\pi$ au sous-groupe $H$.
\ssq
Soit $(\uple X m)$ une base de l'alg\`ebre de Lie $\g h$ du sous-groupe $H$, que l'on compl\`ete en une base $\uple X n$ de l'alg\`ebre de Lie $\g g$ de $G$. Consid\'erant les laplaciens respectifs $\Delta_G=-(X_1^2+\cdots +X_n^2)$ et $\Delta_H=-(X_1^2+\cdots +X_m^2)$ du groupe $G$ et du sous-groupe $H$, l'espace $\Cal H_\pi^s$ (resp. $\Cal H_{\pi\srestr H}^s$) est le domaine de l'op\'erateur $\pi(1+\Delta_G)^s$ (resp. $\pi(1+\Delta_H)^s$). On a donc les inclusions suivantes~:
\prop{II.5.1}
$$\eqalign{& \hbox to 6cm
{$\Cal H_{\pi}^\infty  \subset \Cal H_{\pi\srestr H}^\infty$\hfill}
           \Cal H_{\pi\srestr H}^{-\infty} \subset\Cal H_\pi^{-\infty}\cr             &        \hbox to 6cm
{$\Cal H_{\pi}^s  \subset \Cal H_{\pi\srestr H}^s$ si $s\ge 0$\hfill}  
      \Cal H_{\pi\srestr H}^s   \subset\Cal H_{\pi}^s \hbox{ si }s\le 0
         \cr}$$
\ndem
On peut se demander inversement \`a quelle condition un vecteur $C^\infty$ pour $\pi\restr H$ est un vecteur $C^\infty$ pour $\pi$, et \`a quelle condition un vecteur-distribution pour $\pi$ est un vecteur-distribution pour la restriction $\pi\restr H$~:
\th{II.5.2}
1). Soit $u\in \Cal H_{\pi\srestr H}^\infty$. Consid\'erant alors $u$ comme un vecteur-distribution de $\pi$ on a l'inclusion~:
$$WF(u)\subset \g h^\perp.$$
2). Soit $v\in\Cal H_\pi^{-\infty}$. Alors si $WF(v)\cap \g h^\perp=\emptyset$ le vecteur-distribution $v$ est un vecteur-distribution pour la restriction $\pi\restr H$.
\dem
Soit $\xi\notin \g h^\perp$. Alors il existe un voisinage conique $\Cal V$ de $\xi$ tel que le symbole polynomial $q(\xi)=1+\xi_1^2+\cdots +\xi_m^2$ soit $\Cal V$-elliptique d'ordre $2$, c'est-\`a-dire qu'il v\'erifie~:
$$C_1\Lambda(\xi)^2\le q(\xi)\le C_2\Lambda(\xi)^2$$
pour tout $\xi\in\Cal V$. Pour tout entier positif $k$ le symbole $q^{\#k}$ est alors $\Cal V$-elliptique d'ordre $2k$. Il existe donc un symbole $s_k$ tel que $ES(s_k)\cap \Cal V=\emptyset$ et tel que~:
$$q^{\#k}=r_k+s_k$$
o\`u $r_k$ est elliptique d'ordre $2k$. Soit $r_{-k}$ une param\'etrixe pour $r_k$. Soit $u\in\Cal H_{\pi\srestr H}^\infty$. On a alors~:
$$r_{-k}^{W,\pi}(q^{\#k})^{W,\pi}u = (r_{-k}\# r_k)^{W,\pi}u+(r_{-k}\# s_k)^{W,\pi}u,$$
soit encore~:
$$u=-(r_{-k}\# s_k)^{W,\pi}u + r_{-k}^{W,\pi}(q^{\#k})^{W,\pi}u + (1-r_{-k}\# r_k)^{W,\pi}u.$$
$\xi$ n'appartient pas au front d'onde du premier terme gr\^ace au 2) de la proposition II.3.1, le deuxi\`eme terme appartient \`a $\Cal H_\pi^{2k}$ gr\^ace \`a la proposition II.2.1 et gr\^ace au fait que $(q^{\#k})^{W,\pi}u$ appartient \`a $\Cal H_\pi$ pour tout $k$, et enfin le troisi\`eme terme appartient \`a $\Cal H_\pi^\infty$ par d\'efinition d'une param\'etrixe. On a donc~:
$$u=u_1+u_2$$
o\`u $\xi\notin WF(u_1)$ et $u_2\in \Cal H_{2k},$
ce qui \'equivaut au fait (cf prop. II.3.3) que $\xi$ n'appartient pas \`a $WF_{2k}(u)$. Ceci \'etant vrai pour tout $k$ on en d\'eduit le 1) du th\'eor\`eme.
\ssq
Pour d\'emontrer le 2) on consid\`ere un $v\in\Cal H_\pi^{-\infty}$ tel que $WF(v)\cap \g h^\perp=\emptyset$. D'apr\`es la proposition II.3.1 il existe alors un symbole $p=1+r$ avec $p^{W,\pi}v\in\Cal H_\pi^\infty$ et $ES(r)\cap \g h^\perp=\emptyset$. Pour tout entier $k\ge 0$ on consid\`ere alors une param\'etrixe $\delta_k$ pour le symbole $(1+\xi_1^2+\cdots +\xi_m^2)^{\#k}$ vu comme elliptique sur $\g h^*$, et on en fait un symbole sur $\g g^*$ en posant~:
$$\varepsilon_k(\xi)=\delta_k(\xi+\g h^\perp)$$
moyennant l'identification de $\g h^*$ avec $\g g^*/\g h^\perp$. On voit alors que $\varepsilon_k\#r$ est un symbole d'ordre $-2k$ sur $\g g^*$, et donc d'apr\`es la proposition II.2.1 $(1-\Delta_H)^{-k}r^{W,\pi}$ envoie $\Cal H_\pi^s$ dans $\Cal H_\pi^{s+2k}$. On \'ecrit alors~:
$$\delta_k^{W,\pi}v=\delta_k^{W,\pi}p^{W,\pi}v-(\delta_k\#r)^{W,\pi}v$$
d'o\`u finalement, si $v$ appartient \`a $\Cal H_\pi^s$ on a~:
$$(1-\Delta_H)^{-k}v\in \Cal H_\pi^{s+2k}.$$
d'o\`u $(1-\Delta_H)^{-k}v\in \Cal H_\pi$ pour $k$ assez grand, ce qui veut dire que $v$ est un vecteur-distribution pour la restriction $\pi\restr H$.
\paragraphe{III. Propagation des singularit\'es}
\alinea{III.1. Symboles classiques}
\qquad On d\'esigne par $\mop{Cl}^m(\g g^*)$ la classe des symboles classiques d'ordre $m$, c'est-\`a-dire qui admettent un d\'eveloppement asymptotique~:
$$p\simeq p_m+p_{m-1}+\cdots$$
o\`u $p_{m-j}$ est positivement homog\`ene de degr\'e $m-j$. On note $A\mop{Cl}^{m,Q}$ et $A\mop {Cl}^m$ les intersections $\mop{Cl}^m\cap AS^{m,Q}_1$ et 
 $\mop{Cl}^m\cap AS^{m}_1$ respectivement. D'apr\`es le th\'eor\`eme d'approximation il existe pour tout $p\in\mop{Cl}^m$ un symbole $\widetilde p$ dans $A\mop{Cl}^{m,Q}$ ayant le m\^eme d\'eveloppement asymptotique que $p$. On appelle symbole principal le terme positivement homog\`ene de plus haut degr\'e.
\alinea{III.2. Une forme faible de l'in\'egalit\'e de G\aa rding pr\'ecis\'ee}
\qquad On peut \'etablir une in\'egalit\'e de G\aa rding qui s'\'enonce ainsi~: lorsque $p$ est un symbole hypoelliptique positif il existe une constante $C$ telle que $p^{W,\pi}+CI$ soit un op\'erateur positif \cite {M2 prop. I.4.1}. Lorsque $p$ est seulement suppos\'e positif sans hypoth\`ese d'hypoellipticit\'e un substitut possible est une forme faible de l'{\sl in\'egalit\'e de G\aa rding pr\'ecis\'ee\/}~:
\th{III.2.1}
Soit $p\in AS^{m,Q}_\rho(\g g^*)$ avec $\rho>\frac 12$ et $Q$ assez petit. On suppose $p(\xi)\ge 0$ pour tout $\xi\in\g g^*$. Alors pour tout $\varepsilon>0$ il existe une constante $C_\varepsilon$ telle que pour tout $u\in\Cal H_\pi^\infty$~:
$$<p^{W,\pi}u,u>\ge -C_\varepsilon\|u\|^2_{m-(\rho-1/2-\varepsilon)\over 2}.$$
\dem
On se ram\`ene \`a l'in\'egalit\'e de G\aa rding simple de la fa\c con suivante~: soit $m_0=m-(\rho-\frac 12 -\varepsilon)$. On suppose $\varepsilon\le \rho-\frac 12$, de sorte que $m_0\le m$. On pose~:
$$\widetilde p_{m_0}(\xi)=p(\xi)+\Lambda(\xi)^{m_0}$$
et pose $p_{m_0}=T^Q(\widetilde p_{m_0})$ (cf \S\ I.1). On a l'encadrement~:
$$C_1 \Lambda(\xi)^{m_0}\le p_{m_0}(\xi)\le C_2 \Lambda(\xi)^m,$$
et les estimations~:
$$|D^\alpha p_{m_0}(\xi)|\le C'_\alpha \Lambda(\xi)^{(\rho-m+m_{0})|\alpha|}$$
ce qui veut dire que $p_{m_0}$ est un symbole hypoelliptique appartenant \`a la classe~:
$$AHS_\delta^{m,m_0,Q}(\g g^*)$$
avec $\delta=\rho-m+m_0=\frac 12+\varepsilon$. Il existe donc \cite {M2 prop. I.3.7} une racine carr\'ee approch\'ee $q_{m_0}\in AHS_\delta^{m/2,m_{0}/2,Q}$ telle que~:
$$q_{m_0}\#q_{m_0}=p_{m_0}+r_{m_0}$$
avec $r_{m_0}\in S(\g g^*)$. On en d\'eduit~:
$$<p_{m_0}^{W,\pi}u,\, u>\ge <r_{m_0}^{W,\pi}u,\,u>,$$
d'o\`u finalement~:
$$\eqalign{<p^{W,\pi}u,u>       &=<p_{m_0}^{W,\pi}u,\, u>+<(\Lambda^{m_0})^{W,\pi}u,\, u>\cr
                                &\ge C_{m_0}\|u\|^2_{m_0/2}
        +C'_{m_0,s}\|u\|_s^2\cr}$$
d'apr\`es la proposition II.2.1. Posant $s=m_0/2$ on a donc le r\'esultat.
\qed
\alinea{III.3. Th\'eor\`eme de propagation des singularit\'es}
\qquad On rappelle qu'une {\sl bicaract\'eristique\/} d'une fonction $a\in C^\infty(\g g^*)$ est une courbe int\'egrale du champ hamiltonien $H_a$. Comme on a~:
$$H_a.a=\{a,a\}=0,$$
la fonction $a$ est constante le long de ses bicaract\'eristiques. Une {\sl bicaract\'eristique nulle\/} est une bicaract\'eristique sur laquelle la fonction $a$ s'annule.
\ssq
Nous pouvons maintenant \'enoncer le r\'esultat principal, dont la d\'emonstration occupera le reste du paragraphe~:
\th{III.3.1 \rm(Propagation des singularit\'es)}
Soit $p\in A\mop{Cl}^m(\g g^*)$, soit $a$ la partie r\'eelle du symbole principal $p_m$ de $p$. On se donne un segment $\gamma:[t_0,t_1]\rightarrow \g g^*-\{0\}$ de bicaract\'eristique nulle pour $a$, et on suppose que la partie imaginaire de $p_m$ est positive sur un voisinage de $\gamma$. 
\ssq
Soit $u$ un vecteur-distribution tel que $WF_s(p^{W,\pi}u)\cap \gamma=\emptyset$. Alors pour tout $\delta\in ]0,1/2[$, si $\gamma(t_1)\notin WF_{s+m-\delta}(u)$ on a~: $WF_{s+m-\delta}(u)\cap \gamma=\emptyset$.
\dem
Nous suivons de pr\`es la d\'emonstration de M.E. Taylor \cite {T. \S\ VI.2} pour une distribution sur une vari\'et\'e. La seule diff\'erence notable est la n\'ecessit\'e de se limiter \`a $\delta<1/2$ dans notre contexte, alors que l'on peut prendre $\delta=1$ dans le cas d'une distribution. Cette restriction est d\'ej\`a pr\'esente dans la forme faible de l'in\'egalit\'e de G\aa rding pr\'ecis\'ee (Th\'eor\`eme II.4.1), et provient en d\'efinitive de la restriction ``$\rho>1/2$'' dans le calcul symbolique.
\ssq
On peut se ramener au cas o\`u $m=\delta$~: pour ce faire on multiplie $p^{W,\pi}$ \`a gauche par un op\'erateurs $e^{W,\pi}$ o\`u $e$ est r\'eel elliptique d'ordre $\delta-m$. Le symbole principal de $e\# p$ est $e_{\delta-m}p_m$, sa partie r\'eelle est $ae_m$, et il est facile de voir que les bicaract\'eristiques de $a$ et de $ap_m$ sont les m\^emes.
\ssq
On suppose donc $p$ d'ordre $\delta$. Par ailleurs si le th\'eor\`eme est d\'emontr\'e pour une valeur particuli\`ere de $s$ on l'obtient pour les autres valeurs en appliquant \`a $u$ un op\'erateur $q^{W,\pi}$ avec $q$ elliptique d'ordre bien choisi. Dans la d\'emonstration on peut donc supposer que $s$ est tel que $u\in \Cal H_\pi^{s-\delta /2}$.
\ssq
Soit donc $\gamma:[t_0,t_1]\rightarrow \g g^*-\{0\}$ un segment de bicaract\'eristique nulle pour $a={\Cal R e}\,p_\delta$, soit $\Cal C$ un voisinage conique de $\gamma$, et soit $\Cal M$ une partie born\'ee de $AS^{s}_1(\g g^*)$ constitu\'ee de symboles de $AS^{s-\delta}_1(\g g^*)$  \`a valeurs r\'eelles, d'ordre $s-\delta$, et dont le support essentiel est contenu dans $\Cal C$. On consid\`ere la d\'ecomposition~:
$$p=a+ib$$
avec $b\ge 0$ sur $\Cal C$.
\lemme{III.3.2}
Il existe deux constantes $C_1$ et $C_2$ telles que pour tout $c\in\Cal M$ on a l'encadrement~:
$$\displaylines{{1\over 2}<\{a,c^2\}^{W,\pi}u,\,u>-C_1\|c^{W,\pi}u\|^2
        \le {\Cal I m}<c^{W,\pi}p^{W,\pi}u,\, c^{W,\pi}u>+C_2
                \le C_2+\hfill\cr
        \hfill+{1\over 2}(\|c^{W,\pi}u\|^2
                +\|c^{W,\pi}p^{W,\pi}u\|^2).\cr}$$
\dem
La seconde in\'egalit\'e est imm\'ediate. Pour la premi\`ere on \'ecrit~:
$$\eqalign{{\Cal I m}<c^{W,\pi}p^{W,\pi}u,\, c^{W,\pi}u>
        &={\Cal R e}<b^{W,\pi}c^{W,\pi}u,\, c^{W,\pi}u>\cr
        &+{\Cal R e}< [c^{W,\pi},b^{W,\pi}]u,\, c^{W,\pi}u>\cr
        &+{\Cal I m}<[c^{W,\pi},a^{W,\pi}]u,\, c^{W,\pi}u>,\cr}$$
et on minore successivement les trois termes du membre de droite. Comme $b\ge 0$ sur $\Cal C$, l'application de l'in\'egalit\'e de G\aa rding pr\'ecis\'ee (th\'eor\`eme II.4.1) avec $\varepsilon=1/2-\delta$ donne pour le premier terme~:
$$<b^{W,\pi}c^{W,\pi}u,\, c^{W,\pi}u>\ge -C.\|c^{W,\pi}u\|^2$$
o\`u la constante $C$ ne d\'epend pas de $c\in \Cal M$. Pour la minoration du deuxi\`eme terme on \'ecrit~:
$$<[c^{W,\pi},b^{W\pi }]u,\,c^{W,\pi}u>=<c^{W,\pi}[c^{W,\pi},b^{W\pi }]u,\,u>.$$
Le symbole principal de l'op\'erateur $c^{W,\pi}[c^{W,\pi},b^{W\pi }]$ est 
$ic\{c,b\}$ qui est purement imaginaire et d'ordre $2s-\delta-1$. Compte tenu des propri\'et\'es du calcul symbolique le symbole principal de la partie r\'eelle est donc d'ordre $2s-\delta-2$. L'op\'erateur $c^{W,\pi}[c^{W,\pi},b^{W\pi }]$ envoie donc contin\^ument $\Cal H^{s-\delta/2}_\pi$ dans $\Cal H^{-s+2+\delta/2}$, et d'apr\`es les hypoth\`eses sur $\Cal M$ on a les estimations~:
$$\|c^{W,\pi}[c^{W,\pi},b^{W\pi }]u\|_{-s+2-\frac {3\delta} 2}
        \le C\|u\|_{s-\frac\delta 2}$$
o\`u $C$ la constante $C$ est ind\'ependante de $c\in\Cal M$. Comme $\delta<1/2$ on a a fortiori~:
$$\|c^{W,\pi}[c^{W,\pi},b^{W\pi }]u\|_{-s+\frac \delta 2}
        \le C\|u\|_{s-\frac \delta 2}.$$
On en d\'eduit donc une minoration pour le deuxi\`eme terme~:
$$\Cal R e <[c^{W,\pi},b^{W,\pi}]u,\, c^{W,\pi}u>\ge -C''$$
o\`u $C''$ est ind\'ependante de $c\in\Cal M$. Enfin pour le troisi\`eme terme on remarque que le symbole principal de l'op\'erateur $c^{W,\pi}[c^{W,\pi},a^{W\pi}]$ est $ic\{c,a\}=\frac i2\{c^2,a\}$ qui est d'ordre $2s-\delta-1$, et que l'op\'erateur
$c^{W,\pi}[c^{W,\pi},a^{W\pi}]-ic\{c,a\}^{W,\pi}$ est d'ordre $2s-\delta-2$. Par un raisonnement analogue on obtient la minoration~:
$$\Cal I m<[c^{W,\pi},a^{W,\pi}]u,\, c^{W,\pi}u>
        \ge {1\over 2}<\{a,c^2\}^{W,\pi}u,\,u>-C'''$$
o\`u la constante $C'''$ est ind\'ependante de $c\in \Cal M$. R\'eunissant les trois minorations on en d\'eduit le lemme.
\qed
\lemme{III.3.3}
Soit $e=\{a,c^2\}-(2C_1+1)c^2$ o\`u $C_1$ est la constante donn\'ee par le lemme III.2.2. Alors il existe une constante $C_3$ ind\'ependante de $c\in\Cal M$ telle que~:
$$\Cal R e<e^{W,\pi}u,\,u>\le C_3.$$
\dem
Cela d\'ecoule facilement du fait que l'ensemble $\{c^{W,\pi}p^{W,\pi}u,\,c\in\Cal M \}$ est born\'e dans $\Cal H_\pi^0$.
\qed
Nous allons maintenant construire une famille $\Cal M$ particuli\`ere v\'erifiant les propri\'et\'es demand\'ees. Nous aurons pour cela besoin du lemme suivant~:
\lemme{III.3.4}
Supposons que le champ hamiltonien $H_a$ soit non radial en $\gamma(t_0)$. Soit $\Cal U$ un voisinage conique de $\gamma(t_1)$. Alors il existe trois symboles $c$, $a_0$, $a_1$ tels que~:
\smallskip
1) $c$ est positivement homog\`ene de degr\'e $s$, a son support dans $\Cal C$, et $c(\xi)\le 0$ pour tout $\xi\in\g g^*-\{0\}$.
\smallskip
2) $\{a,c\}\ge 0$ sur $\Cal C\backslash \Cal U$ et $\{a,c\}> 0$ sur 
$\gamma\backslash \Cal U$.
\smallskip
3) $a_0\in S_1^0(\g g^*)$ et $\{a,a_0\}=1$ sur un voisinage de $\gamma$.
\smallskip
4) $a_1\in S^1_1(\g g^*)$ et $\{a,a_1\}$ ne s'annule pas sur $\gamma$.
\dem
Remarquons tour d'abord que si le hamiltonien $H_a$ est radial en $\gamma(t_0)$, alors par homog\'en\'eit\'e de $H_a$ la bicaract\'eristique passant par $\gamma(t_0)$ est elle aussi radiale. Comme les fronts d'onde sont par d\'efinition coniques le th\'eor\`eme III.3.1 est trivial dans ce cas-l\`a.
\ssq
 Supposons donc le champ hamiltonien $H_a$ non radial en $\gamma(t_0)$. Il existe alors une hypersurface conique $\Omega\subset \g g^*-\{0\}$ passant par $\gamma(t_0)$ et transverse \`a \`a $H_a$ dans un voisinage conique de $\gamma(t_0)$. Soit $g\in C^\infty(\g g^*-\{0\})$, positive, positivement homog\`ene de degr\'e $s$, et telle que $g\bigl(\gamma(t_0)\bigr)>0$. On r\'esout alors l'\'equation de transport~:
$$\eqalign{     &H_a.c=g        \cr
                &c\restr \Omega=0.\cr}$$
ce qui nous donne un symbole $c$ positivement homog\`ene de degr\'e $s$ et tel que $\{a,c\}>0$ sur $\gamma$. Quitte \`a le tronquer on obtient facilement un symbole $c$ v\'erifiant les conditions 1 et 2 du lemme. La construction de $a_0$ et $a_1$ se fait de mani\`ere analogue.
\qed
{\sl Fin de la d\'emonstration du th\'eor\`eme\/}~: on d\'efinit la famille $\Cal M$ par~:
$$\Cal M=\{T^Qc_{\lambda,\varepsilon},\,0<\varepsilon\le 1 \}$$
o\`u $\lambda$ est un r\'eel \`a fixer ult\'erieurement et $T^Q$ est l'op\'erateur de r\'egularisation du \S\ I.1. Le symbole $c_{\lambda,\varepsilon}$ est quant \`a lui d\'efini par~:
$$c_{\lambda,\varepsilon}(\xi)=c(\xi)e^{\lambda a_0(\xi)}
        \bigl(1+\varepsilon^2a_1^2(\xi) \bigr)^{-\frac 12}.$$
On pose \'egalement~:
$$\eqalign{e_{\lambda,\varepsilon}      &=\{a,c_{\lambda,\varepsilon}^2\}-(2C_1-1)c_{\lambda,\varepsilon}^2\cr
                                        &=\bigl(\{a,c^2\}
        +(2\lambda-2C_1-1)c^2\bigr)e^{2\lambda a_0}(1+\varepsilon^2a_1^2)\inver.\cr}$$
Fixons $\lambda\ge\frac 12(2C_1+1)$. Le symbole $\bigl(\{a,c^2\}+(\lambda-2C_1-1)c^2\bigr)e^{\lambda a_0}$ est d'ordre $2s$, positif sur $\Cal C\backslash \Cal U$ et strictement positif sur $\gamma\backslash U$. Il existe donc $q$ et $r$ symboles positivement homog\`enes d'ordre $s$ positifs sur $\Cal C$ avec $\mop{supp}q\subset U$ et $r>0$ sur $\Cal C$ tels que~:
$$r^2\le\bigl(\{a,c^2\}+(\lambda-2C_1-1)c^2\bigr)e^{\lambda a_0}+q^2,$$
ce qui s'\'ecrit encore~:
$$q_\varepsilon^2+e_{\lambda,\varepsilon}-r_{\varepsilon}^2\ge 0,$$
avec $r_\varepsilon=r(1+\varepsilon a_1^2)^{-1/2}$ et $q_\varepsilon=q(1+\varepsilon a_1^2)^{-1/2}$.
\ssq
L'in\'egalit\'e de G\aa rding pr\'ecis\'ee (Th\'eor\`eme II.4.1) s'applique au $T^Q(q_\varepsilon^2+e_{\lambda,\varepsilon}-r_\varepsilon^2)$ qui est positif \`a un \'el\'ement de $\Cal S(\g g^*)$ pr\`es. On a donc~:
\let\wt=\widetilde
$$<\wt{r_\varepsilon^2}^{W,\pi}u,\,u>\le <\wt e_{\lambda,\varepsilon}^{W,\pi}u,\,u>+<\wt{q_\varepsilon^2}^{W,\pi}u,\,u>+C_4,$$
(le tilde sur un symbole signifie qu'on lui applique l'op\'erateur de r\'egularisation $T^Q$) d'o\`u on tire facilement~:
$$\eqalign{\|\wt{r_\varepsilon}u\|^2_0  &\le <\wt e_{\lambda,\varepsilon}^{W,\pi}u\,u>
        +\|\wt {q_\varepsilon}^{W,\pi}\|^2_0+C_4        \cr
                                        &\le C_5.\cr}$$
Les constantes $C_4$ et $C_5$ sont ind\'ependantes de $\varepsilon$. En faisant tendre $\varepsilon$ vers z\'ero on a donc~:
$$\|\wt r^{W,\pi}u\|_0^2 \le C_5$$
avec un symbole $\wt r$ qui est $\Cal C$-elliptique, c'est-\`a-dire qui v\'erifie~:
$$C'\Lambda(\xi)^s\le \wt r(\xi)\le C''\Lambda(\xi)^s$$
pour $\|\xi\|$ assez grand. Il existe donc (cf. la d\'emonstration de la proposition II.3.2) un symbole $a$ d'ordre $-s$ tel que~:
$$a\#\wt r=1+b$$
avec $ES(b)\cap \Cal C=\emptyset$. On a donc la d\'ecomposition~:
$$u=a^{W,\pi}r^{W,\pi}u-b^{W,\pi}u$$
avec $a^{W,\pi}r^{W,\pi}u\in \Cal H_\pi^s$ et $\gamma\cap WF(b^{W,\pi}u)=\emptyset$. D'apr\`es la proposition II.3.2 ceci implique que $WF_{s}(u)$ ne rencontre pas $\gamma$.
\qed
Si le symbole principal $p_m$ est r\'eel on peut inverser la direction
du temps. Au vu de l'\'egalit\'e~:
$$WF(u)\backslash WF(p^{W,\pi}u)=\bigcup_t\bigcap_{s\ge t}
        WF_{s+m-\delta}(u)\backslash WF_s(p^{W,\pi}u)$$
on a~:
\cor{III.3.5}
Soit $p\in A\mop{Cl}^m(\g g^*)$ dont la partie principale $p_m$ est
r\'eelle. Alors pour tout vecteur-distribution $u$, pour tout r\'eel
$s$ et pour tout $\delta\in]0,1/2[$, l'ensemble~:
$$WF_{s+m-\delta}(u)\backslash WF_s(p^{W,\pi}u)$$
est invariant par le flot du champ hamiltonien $H_{p_m}$, de m\^eme
que l'ensemble~:
$$WF(u)\backslash WF(p^{W,\pi}u).$$
\ndem
\paragraphe{Appendice~:  calcul symbolique \`a param\`etre}
\qquad Nous d\'emontrons ici le th\'eor\`eme I.2.4 en adaptant la 
d\'emonstration du th\'eor\`eme 4.2 de \cite{M1}. Soit $V$ un espace 
vectoriel r\'eel de dimension finie et $\Cal P$ un c\^one ferm\'e 
dans un espace vectoriel norm\'e quelconque (typiquement, un secteur 
du plan complexe comme au \S\ I.2). La structure de Fr\'echet de 
l'espace $S^m_{\rho,d}(V,\Cal P)$ est donn\'ee par les semi-normes~:
$$N_{\rho,d}^{m,\alpha}(p_\lambda)=\mopl{sup}_\xi
        \Lambda_d(\xi,\lambda)^{-m+\rho|\alpha|}|D^\alpha 
p_\lambda(\xi)|.$$
On consid\'erera \'egalement la famille de semi-normes index\'ee 
par $\N$~:
$$P_{\rho,d}^{m,k}(p_\lambda)=\mopl{sup}_{|\alpha|\le k}
N_{\rho,d}^{m,\alpha}(p_\lambda).$$
Pour tout voisinage compact $Q$ de $0$ dans $V$ le sous-espace 
$AS^{m,Q}_{\rho,d}(V,\Cal P)$ de $S^m_{\rho,d}(V,\Cal P)$ est 
ferm\'e. On d\'esigne par $AS^{-\infty,Q}(V,\Cal P)$ l'intersection lorsque $m$ parcourt $\R$ des $AS^{m,Q}_{\rho,d}(V,\Cal P)$. C'est l'ensemble des $p_\lambda$ tels que $\Cal F\inver p$ a son support inclus dans $Q$, qui appartiennent \`a $S(V)$ pour tout $\lambda$ et tels que les semi-normes~:
$$\mopl{sup}_\xi \Lambda(\xi)^m|D^\alpha p_\lambda(\xi)|$$
soient \`a d\'ecroissance rapide en $\lambda\in\Cal P$.   
Dans la suite on prendra $V=\g g^*$ ou $V=\g g^*\times\g g^*$.
\lemme{A.1 \rm (cf. \cite {M1 prop. 2.5})}
Soit $Q$ un voisinage compact de $0$ dans $V$, et soit $\psi$ une fonction analytique sur $V^*$ dont la s\'erie enti\`ere \`a l'origine converge sur $Q$. Alors l'op\'erateur $\psi(D)=\Cal F\circ \psi\circ \Cal F\inver$ envoie contin\^ument $AS^{m,Q}_{\rho,d}(V,\Cal P)$ dans $AS^{m,Q}_{\rho,d}(V,\Cal P)$.
\ndem
\qquad Soit maintenant $W$ un voisinage \'etoil\'e de $0$ dans l'alg\`ebre de Lie $\g g$ tel que $(\mop{exp}W)^2$ soit l'image diff\'eomorphe par l'exponentielle d'un voisinage $W^2$ de $0$. Pour tout $t\in [0,1]$ et pour tout $x,y$ dans $W$ on pose~:
$$\eqalign{x\mopl ._t y        &=t\inver\mop{Log}(\exp tx.\exp ty) \cr
                                &=x+y+{t\over 2}[x,y]
        +{t^2\over 12}\bigl([x,[x,y]]+[y,[y,x]]\bigr)+\cdots\cr}$$
En appliquant la formule de Taylor en $t=1$ \`a~:
$$p_\lambda\mopl{\#}_t q_\lambda(\xi)=\int\!\!\!\int_{\sg g\times \sg g}
\Cal F\inver p_\lambda(x)\Cal F\inver q_\lambda(y)e^{-i<x\smopl ._t y,\xi>}\,dx\,dy $$
on obtient le d\'eveloppement asymptotique en op\'erateurs bidiff\'erentiels donn\'e dans l'\'enonc\'e du th\'eor\`eme, avec~:
$$C_k(p_\lambda,q_\lambda)(\xi)={1\over k!}\int\!\!\!\int_{\sg g\times\sg g}
\Cal F\inver p_\lambda(x)\Cal F\inver q_\lambda(y){d^k\over dt^k}\restr{t=0}e^{-i<x\smopl ._t y,\xi>}\,dx\,dy,$$
et
$$R_N(p_\lambda,q_\lambda)(\xi)={1\over N!}\int_0^1(1-t)^N\int\!\!\!\int_{\sg g\times\sg g}
\Cal F\inver p_\lambda(x)\Cal F\inver q_\lambda(y) {d^k\over dt^k}e^{-i<x\smopl ._t y,\xi>}\,dx\,dy\,\,dt.$$
On rappelle le r\'esultat suivant~:
\lemme{A.2 \rm (\cite {M1 lemma 4.3})}
Pour tout $k\ge 0$ on a~:
$${d^k\over dt^k}e^{-i<x\smopl ._t y,\xi>}=\psi_k(x,y,\xi,t)e^{-i<x\smopl ._t y,\xi>}$$
o\`u $\psi_k$ est une fonction analytique polynomiale en $\xi$, dont la s\'erie enti\`ere \`a l'origine converge pour tout $t\in[0,1]$ et $x,y\in W$. On a plus pr\'ecis\'ement~:
$$\psi(x,y,\xi,t)=\sum_{r=1}^k \psi_k^r(x,y,\xi,t)$$
o\`u $\psi_k^r$ est polynomiale homog\`ene de degr\'e $r$ en $\xi$. De plus le d\'eveloppement en s\'erie enti\`ere par rapport \`a $t$ s'\'ecrit~:
$$\psi_k^r(x,y,\xi,t)=\sum_{s\ge k+r}\psi_k^{r,s}(x,y,\xi)t^{s-k-r},$$
o\`u $\psi_k^{r,s}$ est polynomiale en les variables $x,y,\xi$, de valuation $\ge r$ par rapport \`a $x$, de valuation $s$ par rapport \`a $(x,y)$, de valuation $\ge r$ par rapport \`a chacune des variables $x$ et $y$, et homog\`ene de degr\'e $r$ par rapport \`a $\xi$.
\ndem
Le lemme A.2 entra\^\i ne imm\'ediatement que les $C_k$ sont des op\'erateurs bi-diff\'erentiels \`a coefficients polynomiaux qui r\'ealisent pour tout voisinage compact $Q$ de $0$ dans $W$ une correspondance bilin\'eaire continue~:
$$AS_{\rho,d}^{m_1,Q}(\g g^*,\Cal P)\times AS_{\rho,d}^{m_2,Q}(\g g^*,\Cal P)
\longrightarrow AS_{\rho,d}^{m_1+m_2-(2\rho-1)k,Q+Q}(\g g^*,\Cal P).$$
La difficult\'e r\'eside donc dans l'estimation du reste $R_N$. On consid\`ere \cite{Me2} pour tout $t\in [0,1]$ l'application~:
$$\eqalign{S_t:W\times W        &\longrightarrow \g g\times \g g \cr
        (x,y)                   &\longrightarrow {1\over 2}
                (x-y+x\mopl ._t y,\,y-x+x\mopl ._t y).\cr}$$
\qquad Il existe un voisinage ouvert $\Omega\subset W$ tel que pour tout $t\in[0,1]$ cette application $S_t$ soit un diff\'eomorphisme sur son image. On notera $K_t$ le jacobien de $S_t$.
\ssq
Soit $Q$ un voisinage compact de $0$ dans $\g g$ contenu dans $\Omega$. On consid\`ere l'op\'erateur $\Cal B_t$ d\'efini sur $AS^{-\infty, Q\times Q}(\g g^*\times\g g^*,\Cal P)$ par~:
$$\Cal B_t p_\lambda(\xi,\eta)       
        =\int\!\!\!\int_{\sg g\times \sg g}
\Cal F\inver ( p_\lambda)(x,y)e^{-i<S_t(x,y),\,(\xi,\eta)>}\,dx\,dy,$$
de sorte que l'on a~:
$$p_\lambda\mopl{\#}_t q_\lambda(\xi)=\Cal B_t(p_\lambda\otimes q_\lambda)(\xi,\xi).$$
Le reste $R_N$ s'\'ecrit donc, gr\^ace au lemme A.2~:
$$\eqalign{\hbox to -6.1mm{}R_N(p_\lambda,q_\lambda)(\xi)        &={1\over N!}\int_0^1(1-t)^N\!\!
        \int\!\!\!\int_{\sg g\times\sg g}\hbox to -15pt{}\Cal F\inver(p_\lambda\otimes q_\lambda)(x,y)
        \psi_{N+1}(x,y,\xi,t)e^{-i<S_t(x,y),\,(\xi,\xi)>}dxdy\,dt \cr
                                &={1\over N!}\int_0^1(1-t)^N
        \bigl(\Cal B_t\circ\widetilde \psi_{N+1}(D_1,D_2,\xi_1,\xi_2,t)
        \bigr)(p_\lambda\otimes q_\lambda)(\xi,\xi),\cr}$$
o\`u $\widetilde \psi_{N+1}(D_1,D_2,\xi,\eta,t)$ est l'op\'erateur pseudo-diff\'erentiel (au sens classique) sur $\g g^*\times \g g^*$ de symbole \`a gauche~:
$$\widetilde \psi_{N+1}(x,y,\xi,\eta,t)=\psi_{N+1}(x,y,{\xi+\eta\over 2},t).$$
\lemme{A.3 \rm(cf. \cite{M1 prop. 3.6})}
Pour tout $ p_\lambda\in AS^{-\infty, Q\times Q}(\g g^*\times \g g^*,\Cal P)$ on a les estimations~:
$$P^{0,0}_{\rho,d}(\Cal B_t p_\lambda)\le
        C.P_{\rho,d}^{-2n-1,0}( p_\lambda).$$
\dem
On a la suite d'in\'egalit\'es~:
$$\eqalign{|\Cal B_t  p_\lambda (\xi)|
        &\le\int_Q|\Cal F\inver  p_\lambda (x)|\,dx       \cr
        &\le\mop{Vol}Q.\mopl{sup}_x |\Cal F\inver  p_\lambda (x)|\cr
        &\le\mop{Vol}Q.\int_{\sg g^*\times\sg g^*}| p_\lambda(\zeta)|\,d\zeta        \cr
        &\le\mop{Vol}Q. \mopl{sup}_\eta \Lambda_d(\eta,\lambda)^{2n+1}
        | p_\lambda(\eta)|\int_{\sg g^*\times\sg g^*}
        \Lambda_d(\zeta,\lambda)^{-2n-1}\,d\zeta, \cr}$$
d'o\`u le r\'esultat.
\qed
A partir de la d\'efinition de l'op\'erateur $\Cal B_t$ il est facile d'obtenir les deux lemmes suivants~:
\lemme{A.4 \rm [M1 lemma 3.7]}
Pour tout multi-indice $\alpha$ et pour tout $t\in[0,1]$ on a l'\'egalit\'e suivante entre op\'erateurs sur $AS^{-\infty,Q\times Q}(\g g^*\times \g g^*,\Cal P)$~:
$$D^\alpha\circ \Cal B_t=\Cal B_t\circ S_t(D)^\alpha$$
\ndem
\lemme{A.5}
Soit $\xi_k, k=1,\ldots, 2n$ une coordonn\'ee sur $\g g^*\times \g g^*$. Alors~:
$$\Cal B_t\circ \xi_k=\Bigl(\sum_{i=1}^{2n}R_{ki,t}(D)\xi_i\Bigr)\circ\Cal B_t,$$
o\`u $R_{ki,t}(x)$ d\'esigne le coefficient $(k,i)$ de la matrice jacobienne de $S_t$ en $x$.
\ndem 
A partir de l'\'ecriture~: 
$$\Cal B_t p_\lambda(\xi,\eta)       
       =\int\!\!\!\int_{\sg g\times \sg g}
\Cal F\inver ( p_\lambda)\circ S_t\inver(x,y)e^{-i<(x,y),\,(\xi,\eta)>}K_t(x,y)\inver\,dx\,dy$$
obtenue \`a partir de la d\'efinition en effectuant le changement de variable $x'=S_t(x)$ on obtient \'egalement~:
\lemme{A.6 \rm \cite {M1 lemma 3.8}}
Soit $\xi_k, k=1,\ldots, 2n$ une coordonn\'ee sur $\g g^*\times \g g^*$. Alors~:
$$\xi_k\circ \Cal B_t=\Cal B_t\circ\Bigl(
\sum_{i=1}^{2n}(P_{ki,t}\circ S_t)(D)\xi_i+\bigl({D_k K_t\inver
        \over K_t\inver}\circ S_t\bigr)(D)\Bigr),$$
o\`u $P_{ki,t}(x)$ d\'esigne le coefficient $(k,i)$ de la matrice jacobienne de $S_t\inver$ en $x$.
\ndem
En combinant les lemmes A.4 et A.6 avec le lemme A.1 on arrive \`a permuter l'op\'erateur $\Cal B_t$ avec n'importe quel op\'erateur diff\'erentiel \`a coefficients polynomiaux~:
\lemme{A.7 \rm \cite{M1 prop. 3.9}}
Soit $P$ un polyn\^ome \`a $4n$ variables, et soit $P(\xi,D)$ un op\'erateur diff\'erentiel \`a coefficients polynomiaux sur $\g g^*\times \g g^*$ qui admet $P$ comme symbole pour un certain choix d'ordre des variables $(\uple \xi{2n},\uple D{2n})$. Alors si $P$ est de degr\'e $j$ par rapport aux variables $\xi$ et de valuation $v$ par rapport aux variables $D$, pour tout $t\in [0,1]$ on a~:
$$P(\xi,D)\circ \Cal B_t=\Cal B_t\circ P_t$$
o\`u $P_t$ est un op\'erateur (en g\'en\'eral non-diff\'erentiel) qui envoie contin\^ument $AS^{m,Q\times Q}_\rho (\g g^*\times \g g^*,\Cal P)$ dans $AS^{m+j-\rho v,Q\times Q}_\rho (\g g^*\times \g g^*,\Cal P)$.
\ndem\qquad
A partir des lemmes A.3, A.5 et A.7 on obtient les estimations fondamentales pour l'op\'erateur $\Cal B_t$~: 
\lemme{A.8}
Pour tout entier positif $k$ et pour tout r\'eel $m$ il existe un entier $M(k)$ et une constante $C$ tels que pour tout $t\in[0,1]$ et $p\in AS^Q_d(\g g^*\times \g g^*)$ on a~:
$$P_{\rho,d}^{m,k}(\Cal B_tp_\lambda)\le C.P_{\rho,d}^{m,M(k)}(p_\lambda).$$
\dem
On suppose d'abord que $m\le 0$. On consid\`ere pour tout multi-indice $\alpha$ le plus petit entier pair $\mu_\alpha$ plus grand que $-m+\rho|\alpha|$. En appliquant le lemme A.6 aux op\'erateurs~:
$$\Lambda_d(\xi)^{\mu_\alpha}D^\alpha,\hbox to 6mm{}|\alpha|\le k$$  
on obtient les estimations~:
$$P_{\rho,d}^{m,k}(\Cal B_tp_\lambda)\le C.P_{\rho,d}^{m-2n-3,M(k)}(p_\lambda).$$
Si $m>0$ on applique ceci \`a $p'_\lambda(\xi)=\Lambda_d(\xi,\lambda)^{-\mu}p_\lambda(\xi)$ o\`u $\mu$ est le plus petit entier pair plus grand que $m$, et on applique le lemme A.5. On obtient donc les estimations suivantes pour tout r\'eel $m$~:
$$P_{\rho,d}^{m,k}(\Cal B_tp_\lambda)\le C.P_{\rho,d}^{m-2n-5,M(k)}(p_\lambda).$$
Enfin la formule de Taylor \`a l'ordre $N$ appliqu\'ee \`a $t\rightarrow \Cal B_tp_\lambda(\xi)$ permet de d\'evelopper $\Cal B_t$ en une somme d'op\'erateurs diff\'erentiels plus un reste, ce qui permet pour $N$ assez grand de se d\'ebarrasser du d\'ecalage en $2n+5$.
\qed
\prop{A.9}
$\Cal B_t$ se prolonge en un op\'erateur continu de $AS_{\rho,d}^{m,Q\times Q}(\g g^*\times\g g^*,\Cal P)$ dans $AS_{\rho,d}^{m,Q\times Q}(\g g^*\times\g g^*,\Cal P)$.
\dem
Si $AS^{-\infty,Q\times Q}(\g g^*\times\g g^*,\Cal P)$ \'etait dense dans $AS_{\rho,d}^{m,Q\times Q}(\g g^*\times\g g^*,\Cal P)$ on pourrait conclure imm\'ediatement, mais ce n'est pas le cas. On introduit la {\sl topologie de H\"ormander} \cite{Hr1 \S\ 3}, d\'efinie de la fa\c con suivante~: une suite $p_{n,\lambda}$ converge vers $p_\lambda$ au sens de H\"ormander si elle converge vers $p_\lambda$ au sens $C^\infty$ (sans supposer d'uniformit\'e par rapport au param\`etre $\lambda$) en restant born\'ee dans $AS_{\rho,d}^{m,Q\times Q}(\g g^*\times\g g^*,\Cal P)$. On montre facilement que $\Cal B_t$ est continu aussi pour cette topologie, et que $AS^{-\infty,Q\times Q}(\g g^*\times\g g^*,\Cal P)$ est dense dans $AS_{\rho,d}^{m,Q\times Q}(\g g^*\times\g g^*,\Cal P)$ en ce sens, ce qui permet de conclure.
\qed
Le lemme ci-dessous est une cons\'equence du lemme 1, et le suivant d\'ecoule essentiellement de la r\`egle de Leibniz~:
\lemme{A.10}
L'op\'erateur pseudo-diff\'erentiel $\wt\psi_{N+1}(D_1,D_2,\xi,\eta,t)$ envoie
contin\^ument $AS_{\rho,d}^{m,Q\times Q}(\g g^*\times \g g^*,\Cal P)$ dans $AS_{\rho,d}^{m-(N+1)(2\rho-1),Q\times Q}(\g g^*\times \g g^*,\Cal P)$.
\ndem
\lemme{A.11 \rm\cite{M1 prop. 2.10}}
La restriction \`a la diagonale est continue~:
$$AS_{\rho,d}^{m,Q\times Q}(\g g^*\times \g g^*,\Cal P)\longrightarrow AS_{\rho,d}^{m,Q}(\g g^*,\Cal P).$$
\ndem
{\sl Fin de la d\'emonstration du th\'eor\`eme I.2.4\/}~:
Si $(p_\lambda,q_\lambda)\mapsto p_\lambda\otimes q_\lambda$ r\'ealisait une correspondance bilin\'eaire continue~:
$$AS_{\rho,d}^{m_1,Q}(\g g^*,\Cal P)\times AS_{\rho,d}^{m_2,Q}(\g g^*,\Cal P) \longrightarrow AS_{\rho,d}^{m_1+m_2,Q\times Q}(\g g^*\times\g g^*,\Cal P)$$
le th\'eor\`eme I.2.4 serait d\'emontr\'e. L'assertion est malheureusement fausse, mais on a facilement les estimations suivantes pour $p_\lambda,q_\lambda$ dans $AS^{-\infty,Q}(\g g^*,\Cal P)$~:
$$P^{m_1+m_2,k}_{\rho,d}(p_\lambda\otimes q_\lambda)\le
P^{m_1-\rho k,k}_{\rho,d}(p_\lambda).P^{m_2-\rho k,k}_{\rho,d}(q_\lambda).$$
En \'ecrivant le d\'eveloppement asymptotique de $R_N$ en op\'erateurs bidiff\'erentiels~:
$$R_N=\sum_{l=N+1}^{N'}C_l+R_{N'}$$
pour $N'$ assez grand (d\'ependant de $k$) et en appliquant la proposition A.9, le lemme A.10 et le lemme A.11, on voit qu'il existe un entier $M(k)$ tel que~:
$$P^{m_1+m_2,k}_{\rho,d}\bigl(R_N(p_\lambda,q_\lambda)\bigr)
\le C.P^{m_1,M(k)}_{\rho,d}(p_\lambda).P^{m_2,M(k)}_{\rho,d}(q_\lambda).$$
Un argument de densit\'e analogue \`a celui de la d\'emonstration de la proposition A.9 permet alors de conclure.
\qed
\paragraphe{R\'ef\'erences}
\bib{Dui}J.J. Duistermaat, {\sl Fourier integral operators\/}, Courant Institute, 1973 (R\'e\'ed. Progress in Math. 130, Birkh\"auser 1995).
\bib{D-H}J.J. Duistermaat, L. H\"ormander, {\sl Fourier integral operators II\/}, Acta Math. 128 (1972), 183-269. 
\bib{Go}R. Goodman, {\sl Elliptic and subelliptic estimates for operators in an enveloping algebra\/}, Duke Math. J. 47, No 4 (1980), 819-833.
\bib{He}B. Helffer, {\sl Th\'eorie spectrale pour des op\'erateurs globalement elliptiques\/}, Ast\'erisque, Soc. Math. France 1984.
\bib{Hr1}L. H\"ormander, {\sl The Weyl calculus of pseudodifferential operators,\/} Comm. Pure Appl. Math. 32 (1979), 359-443.
\bib{Hr2}L. H\"ormander, {\sl The analysis of linear partial differential operators III}, Springer 1985.
\bib{Hw}R. Howe, {\sl Wave front sets of representations of Lie groups\/}, in {\sl Automorphic forms, representation theory and Arithmetic\/}, Tata Inst. fund. res. st. math. 10, Bombay 1981.
\bib{J\o}P.E.T. J\o rgensen, {\sl Distribution representations of Lie groups\/}, J. Math. Anal. Appl. 65 (1978), 1-19.
\bib{M1}D. Manchon, {\sl Weyl symbolic calculus on any Lie group\/}, Acta Appl. Math. 30 (1993), 159-186.
\bib{M2}D. Manchon, {\sl op\'erateurs pseudodiff\'erentiels et repr\'esentations unitaires des groupes de Lie\/}, Bull. Soc. Math. France 123 (1995), 117-138.
\bib{M3}D. Manchon, {\sl Formule de Weyl pour les groupes de Lie nilpotents}, J. f.d. Reine u. Angew. Math. 418 (1991), 77-129
\bib{M4}D. Manchon, {\sl Distributions \`a support compact et repr\'esentations unitaires}, pr\'epubl. math. RT/9809005, \`a para\^\i tre au J. Lie Theory.
\bib{Me1}A. Melin, {\sl A remark on invariant pseudo-differential operators}, Math. Scand. 30 (1972), 290-296. 
\bib{Me2}A. Melin, {\sl Parametrix constructions for right invariant differential operators on nilpotent groups}, Ann. Glob. Anal. Geom. 1, N.1 (1983), 79-130.
\bib{Ne}E. Nelson, {\sl Analytic vectors\/}, Ann. Math. 70 No 3 (1959), 572-615.
\bib{S}R.T. Seeley, {\sl Complex powers of an elliptic operator\/}, in {\sl Singular integrals\/}, proc. symp. pure math. 10 (1967), 288-307.
\bib{Shu}M.A. Shubin, {\sl Pseudodifferential operators and spectral theory\/}, Springer 1987.
\bib{St}R.S. Strichartz, {\sl A functional calculus for elliptic pseudo-differential operators\/}, Amer. J. Math. 94 (1972), 711-722.
\bib{Stk}H. Stetk\ae r, {\sl Invariant pseudo-differential operators}, Math. Scand. 28 (1971), 105-123.
\bib{T}M.E. Taylor, {\sl Pseudodifferential operators\/}, Princeton 1981.
\bye